\documentclass[12pt]{amsart}

\usepackage{enumitem}
\usepackage{psfrag}
\usepackage{graphicx}
\usepackage{pinlabel}
\usepackage{graphicx}
\usepackage{amsmath}
\usepackage{amssymb}
\usepackage{amscd}
\usepackage{autobreak}
\usepackage{picinpar}
\usepackage{times}
\usepackage{pb-diagram}
\usepackage{graphicx}
\usepackage{wrapfig}
\usepackage{xspace}
\usepackage{hyperref}
\usepackage{url}
\usepackage{caption}
\usepackage{subcaption}
\usepackage{fancyhdr}
\usepackage{overpic}
\usepackage{color}
\usepackage[square,comma,sort&compress,numbers]{natbib}
\usepackage{latexsym}
\usepackage{mathrsfs}
\usepackage{amsfonts}
\usepackage{hyperref}
\usepackage{amsthm}
\usepackage{appendix}
\usepackage{pdfsync}

\theoremstyle{definition}
\newtheorem{theorem}{Theorem}[section]
\newtheorem{lemma}[theorem]{Lemma}
\newtheorem{proposition}[theorem]{Proposition}

\newtheorem{definition}[theorem]{Definition}

\newtheorem{remark}[theorem]{Remark}
\newtheorem*{theorem*}{Theorem}
\graphicspath{{figures/}}   

\def\qed{\hfill{Q.E.D.}\smallskip}

\begin{document}

\title{\bf Discrete conformal structures on surfaces with boundary (I)---Classification}
\author{Xu Xu, Chao Zheng}

\date{\today}

\address{School of Mathematics and Statistics, Wuhan University, Wuhan, 430072, P.R.China}
 \email{xuxu2@whu.edu.cn}

\address{School of Mathematics and Statistics, Wuhan University, Wuhan 430072, P.R. China}
\email{czheng@whu.edu.cn}

\thanks{MSC (2020): 52C25,52C26}

\keywords{Discrete conformal structures; Classification; Surfaces with boundary; Polyhedral metrics}

\begin{abstract}
In this paper, we introduce the discrete conformal structures on surfaces with boundary in an axiomatic approach, which ensures that the Poincar\'{e} dual of an ideally triangulated surface with boundary has a good geometric structure.
Then we classify the discrete conformal structures on surfaces with boundary,
which turns out to unify and generalize Guo-Luo's generalized circle packings \cite{GL2}, Guo's vertex scalings \cite{Guo} and Xu's partial discrete conformal structures \cite{Xu22} on surfaces with boundary.
This generalizes the results of Glickenstein-Thomas \cite{GT} on closed surfaces to surfaces with boundary.
Motivated by \cite{BPS, Zhang-Guo-Zeng-Luo-Yau-Gu}, we further study the relationships between the discrete conformal structures on surfaces with boundary and the hyperbolic trigonometry.
Unexpectedly, we find that some subclasses of the discrete conformal structures on surfaces with boundary are closely related to the twisted generalized hyperbolic triangles introduced by Roger-Yang \cite{Roger-Yang}, which does not appear in the case of closed surfaces.
Finally, we study the relationships between the discrete conformal structures on surfaces with boundary and the 3-dimensional hyperbolic geometry by constructing ten types of generalized hyperbolic tetrahedra. Some new generalized hyperbolic tetrahedra recently introduced by
Belletti-Yang \cite{B-Y} naturally appears in the constructions.
\end{abstract}

\maketitle

\section{Introduction}

Discrete conformal structures on polyhedral surfaces are discrete analogues of the smooth conformal structures on Riemann surfaces, which assign the discrete metrics defined on the edges by scalar functions defined on the vertices.
In the last two decades,
there have been lots of researches on the discrete conformal structures on closed surfaces.
Different from the case on closed surfaces,
there are few researches on the discrete conformal structures on surfaces with boundary (bordered surfaces sometimes).
The known discrete conformal structures on surfaces with boundary include
Guo-Luo's generalized circle packings \cite{GL2}, Guo's vertex scalings \cite{Guo} and Xu's discrete conformal structures \cite{Xu22}.
Motivated by Thurston's circle packings on closed surfaces \cite{Thurston}, Guo-Luo \cite{GL2} introduced some generalized circle packings on surfaces with boundary.
Following Luo's vertex scalings of piecewise linear metrics on closed surfaces \cite{Luo1},
Guo \cite{Guo} introduced a class of discrete conformal structures, also called vertex scalings, on surfaces with boundary.
Motivated by the relationships between the discrete conformal structures on closed surfaces and the 3-dimensional hyperbolic geometry \cite{BPS, Zhang-Guo-Zeng-Luo-Yau-Gu},
Xu \cite{Xu22} introduced a new class of discrete conformal structures on surfaces with boundary.
A natural question is whether there exist other types of discrete conformal structures on surfaces with boundary and whether there exists
a classification of the discrete conformal structures on surfaces with boundary.
In this paper, we introduce the discrete conformal structures on ideally triangulated surfaces with boundary in an axiomatic approach
pioneered by Glickenstein \cite{Glickenstein}, and
give such a classification of the discrete conformal structures on surfaces with boundary.
In the following works \cite{X-Z DCS2,X-Z DCS3}, we study the geometry of these discrete conformal structures on surfaces with boundary.

\subsection{Definition of discrete conformal structures on surfaces with boundary}
Suppose $\Sigma$ is a compact surface with the boundary $\partial\Sigma$ consisting of $N$ connected components, which are topologically circles.
Let $\widetilde{\Sigma}$ be the compact surface obtained by coning off each boundary component of $\Sigma$ to be a point.
Then there are exactly $N$ cone points $\{v_1,...,v_N\}$ in $\widetilde{\Sigma}$ so that $\widetilde{\Sigma}-\{v_1,...,v_N\}$ is homeomorphic to $\Sigma-\partial\Sigma$.
An ideal triangulation $\mathcal{T}$ of $\Sigma$ is a triangulation $\widetilde{\mathcal{T}}$ of $\widetilde{\Sigma}$ such that the vertices of the triangulation $\widetilde{\mathcal{T}}$ are exactly the cone points $\{v_1,...,v_N\}$.
The ideal edges and ideal faces of $\Sigma$ in the ideal triangulation $\mathcal{T}$ are defined to be the intersections $\widetilde{E}\cap \Sigma$ and $\widetilde{F}\cap \Sigma$ respectively, where $\widetilde{E}$ and $\widetilde{F}$ are the sets of unoriented edges and faces in the triangulation $\widetilde{\mathcal{T}}$ of $\widetilde{\Sigma}$.
The intersection of an ideal face and the boundary $\partial\Sigma$ is called a boundary arc.
For simplicity, we use $B=\{1, 2,...,N\}$ to denote the boundary components and label them as $i\in B$.
We use $E$ to denote the set of unoriented ideal edges and label an unoriented ideal edge between two adjacent boundary components $i,j\in B$ as $\{ij\}$.
We use $E_{+}$ to denote the set of oriented edges and label them as ordered pairs $(i,j)$.
We use $F$ to denote the set of ideal faces and label the ideal face adjacent to the boundary components $i,j,k\in B$ as $\{ijk\}$.
In fact, the ideal face $\{ijk\}$ is a hexagon, which is in one-to-one correspondence with the triangle $v_iv_jv_k$ in $\widetilde{\mathcal{T}}$.
The sets of real valued functions on $B, E$ and $E_{+}$ are denoted by $\mathbb{R}^N$, $\mathbb{R}^E$ and $\mathbb{R}^{E_+}$ respectively.

The edge length function associated to $\mathcal{T}$ is a vector $l\in \mathbb{R}^E_{>0}$ assigning each ideal edge $\{ij\}\in E$ a positive number $l_{ij}=l_{ji}$.
A well-known fact in hyperbolic geometry is that for any three positive numbers, there exists a unique right-angled hyperbolic hexagon (hyper-ideal hyperbolic triangle) up to isometry with the lengths of three non-pairwise adjacent edges given by the three positive numbers (\cite{Ratcliffe} Theorem 3.5.14).
Hence, for an ideal face $\{ijk\}\in F$, there exists a unique right-angled hyperbolic hexagon whose three non-pairwise adjacent edges have lengths $l_{ij}, l_{ik}, l_{jk}$.
Gluing all such right-angled hyperbolic hexagons isomorphically along the ideal edges in pairs, one can construct a hyperbolic surface with totally geodesic boundary from the ideal triangulation $\mathcal{T}$.
Conversely, any ideally triangulated hyperbolic surface $(\Sigma, \mathcal{T})$ with totally geodesic boundary produces a function $l\in \mathbb{R}^E_{>0}$ with $l_{ij}$ given by the length of the shortest geodesic connecting the boundary components $i,j\in B$.
The edge length function $l\in \mathbb{R}^E_{>0}$ is called a \textit{discrete hyperbolic metric} on $(\Sigma,\mathcal{T})$.

Motivated by the partial edge length introduced by Glickenstein \cite{Glickenstein} and Glickenstein-Thomas \cite{GT} for triangulated closed surfaces,
we introduce the following definition of partial edge length for ideally triangulated surfaces with boundary $(\Sigma, \mathcal{T})$.

\begin{definition}\label{Def: partial edge length}
Let $(\Sigma,\mathcal{T})$ be an ideally triangulated surface with boundary.
An assignment of partial edge lengths is a map $d\in \mathbb{R}^{E_{+}}$ such that
$l_{ij}=d_{ij}+d_{ji}>0$ for every edge $\{ij\}\in E$ and
\begin{equation}\label{Eq: compatible condition}
\sinh d_{ij} \sinh d_{jk} \sinh d_{ki}=
\sinh d_{ji} \sinh d_{kj} \sinh d_{ik}
\end{equation}
for every ideal face $\{ijk\}\in F$.
\end{definition}

\begin{figure}[!ht]
  \centering
  \includegraphics[scale=0.9]{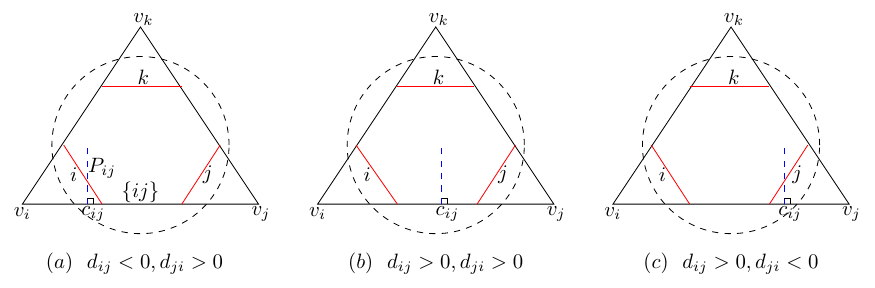}
  \caption{Poincar\'{e} dual to the ideal edge $\{ij\}\in E$ with the position determined by $d_{ij}$ and $d_{ji}$ in the Klein model (dotted).
The red lines represent the boundary arcs of the ideal face $\{ijk\}$.
The black triangle $v_iv_jv_k$ is the corresponding  hyper-ideal hyperbolic triangle.}
\label{Figure 1}
\end{figure}

In general, the partial edge lengths do not satisfy the symmetry condition, i.e., $d_{ij}\neq d_{ji}$.
The condition (\ref{Eq: compatible condition}) in Definition \ref{Def: partial edge length} ensures the existence of a dual geometric structure on the Poincar\'{e} dual of an ideally triangulated surface with boundary $(\Sigma,\mathcal{T})$.
Denote $E_{ij}$ as the hyperbolic geodesic line in the hyperbolic plane $\mathbb{H}^2$ extending the ideal edge $\{ij\}\in E$.
The edge center $c_{ij}$  of the ideal edge $\{ij\}$ is the unique point in $E_{ij}$ that is of signed distance $d_{ij}$ to $i\in B$ and $d_{ji}$ to $j\in B$.
Note that the signed distance $d_{ij}$ of $c_{ij}$ to  $i$ is positive if $c_{ij}$ is on the same side of $i$ as $j$ along the hyperbolic geodesic line $E_{ij}$, and negative otherwise.
Please refer to Figure \ref{Figure 1}.

The hyperbolic geodesic line orthogonal to the hyperbolic geodesic line $E_{ij}$ at $c_{ij}$ is called the edge perpendicular, denoted by $P_{ij}$.
We can identify the perpendicular $P_{ij}$ with a unique 2-dimensional time-like vector subspace $\widetilde{P}_{ij}$ of the Lorentzian space $\mathbb{R}^3$ such that $\widetilde{P}_{ij}\cap \mathbb{H}^2=P_{ij}$.
Since the 2-dimensional subspaces $\widetilde{P}_{ij}$ and $\widetilde{P}_{jk}$ associated to $P_{ij}$ and $P_{jk}$ intersect in a 1-dimensional subspace of $\mathbb{R}^3$,
then the intersection $\widetilde{P}_{ij}\cap\widetilde{P}_{jk}$ corresponds to a point $q$ in the Klein model.
We view this point $q$ as the intersection point of the perpendiculars $P_{ij}$ and $P_{jk}$.
Note that the point $q$ may be not in $\mathbb{H}^2$, i.e., the perpendiculars $P_{ij}, P_{jk}$ may not intersect within $\mathbb{H}^2$.
Furthermore, if the perpendiculars $P_{ij}, P_{jk}, P_{ki}$ of a right-angled hyperbolic hexagon $\{ijk\}$ intersect at a common point $c_{ijk}$, then such an assignment of partial edge lengths naturally induces a geometric structure on the Poincar\'{e} dual of the ideal face $\{ijk\}$.
The point $c_{ijk}$ is called the \textit{geometric center} of the ideal face $\{ijk\}$ (face center for simplicity).
Similarly, $c_{ijk}$ may not lie in $\mathbb{H}^2$ even if it exists.
Lemma \ref{Lem: compatible condition} shows that the perpendiculars $P_{ij}, P_{jk}, P_{ki}$ of $\{ijk\}$ intersect at a common point $c_{ijk}$ if and only if the condition (\ref{Eq: compatible condition}) holds.
This is the motivation of Definition \ref{Def: partial edge length} for the partial edge length on surfaces with boundary.
Please refer to \cite{G_05, Glickenstein, G_24, GT, Thomas} for more information on partial edge lengths and geometric centers.

A discrete conformal structure on an ideally triangulated surface with boundary assigns the partial edge lengths defined on the oriented edges $E_+$ by scalar functions defined on the boundary components $B=\{1, 2,...,N\}$.
Motivated by Glickenstein's axiomatic approach to the discrete conformal structures on triangulated closed surfaces \cite{Glickenstein, GT}, we introduce the following definition of discrete conformal structures on ideally triangulated surfaces with boundary.

\begin{definition}\label{Def: DCS}
Let $(\Sigma,\mathcal{T})$ be an ideally triangulated surface with boundary.
A discrete conformal structure $d=d(f)$ on $(\Sigma,\mathcal{T})$ is a smooth map, sending a function $f\in \mathbb{R}^N$ defined on boundary components $B$ to a partial edge length function $d\in \mathbb{R}^{E_+}$, such that
\begin{equation}\label{Eq: variation 1}
\frac{\partial l_{ij}}{\partial f_i}
=\coth d_{ij}
\end{equation}
for each $(i,j)\in E_{+}$ and
\begin{equation}\label{Eq: variation 2}
\frac{\partial d_{ij}}{\partial f_k}=0
\end{equation}
if $k\neq i$ and $k\neq j$.
\end{definition}

\begin{remark}
The geometric motivation for Definition \ref{Def: DCS} of the discrete conformal structures on surfaces with boundary is as follows.
For a hyper-ideal hyperbolic triangle $v_iv_jv_k$,
if we keep $f_i$ and $f_j$ fixed and perturb $f_k$ in a way satisfying (\ref{Eq: variation 1}) and (\ref{Eq: variation 2}) so that
the vertex $v_k$ is moved to $\overline{v}_k$ infinitesimally,
then the vertex $v_k$, the original face center $c_{ijk}$ and the vertex $\overline{v}_k$ of the new hyper-ideal hyperbolic triangle lie in a 2-dimensional vector subspace of the Lorentzian space  $\mathbb{R}^3$.
In the Klein model, this is equivalent to $v_k$, $c_{ijk}$ and $\overline{v}_k$ stays in a line.
Please refer to Figure \ref{Figure_5}.
Such a variation is called an $f_k-$conformal variation by Glickenstein-Thomas \cite{GT}.
Please refer to Subsection \ref{subsection 2} for more details.
\end{remark}

\subsection{Classification of discrete conformal structures on surfaces with boundary}
A natural question on the discrete conformal structures in Definition \ref{Def: DCS} is
whether the discrete conformal structures contain the known discrete conformal structures in
\cite{GL2, Guo, Xu22}
and whether the discrete conformal structures can be classified.
The following theorem answers this question affirmatively and further gives the explicit forms of the discrete conformal structures in Definition \ref{Def: DCS}.

\begin{theorem}\label{Thm: DCS}
Let $(\Sigma,\mathcal{T})$ be an ideally triangulated surface with boundary and $d=d(f)$ be a discrete conformal structure on $(\Sigma,\mathcal{T})$.
There exist constant vectors $\alpha\in \mathbb{R}^N$ and $\eta\in \mathbb{R}^{E}$ satisfying $\eta_{ij}=\eta_{ji}$, such that for any right-angled hyperbolic hexagon $\{ijk\}\in F$,
\begin{description}
  \item[(A)] if $\frac{\sinh d_{ij}}{\sinh d_{ji}}>0$, then the discrete conformal structure $d=d(f)$ has one of the following two forms
\begin{description}
  \item[(i)]
\begin{equation*}\label{Eq: d3}
\coth d_{ij}
=-\frac{\alpha_ie^{2f_i}}{\sinh l_{ij}}\sqrt{\frac{1+\alpha_je^{2f_j}}{1+\alpha_ie^{2f_i}}}
+\frac{\eta_{ij}e^{f_i+f_j}}{\sinh l_{ij}}
\end{equation*}
with
\begin{equation}\label{Eq: DCS3}
\cosh l_{ij}
=-\sqrt{(1+\alpha_ie^{2f_i})(1+\alpha_je^{2f_j})}
+\eta_{ij}e^{f_i+f_j},\ \text{for}\ 1+\alpha_ie^{2f_i}>0, 1+\alpha_je^{2f_j}>0,
\end{equation}
\begin{equation}\label{Eq: new 1}
\cosh l_{ij}
=\sqrt{(1+\alpha_ie^{2f_i})(1+\alpha_je^{2f_j})}
+\eta_{ij}e^{f_i+f_j},\ \text{for}\ 1+\alpha_ie^{2f_i}<0, 1+\alpha_je^{2f_j}<0;
\end{equation}
\item[(ii)]
\begin{equation*}\label{Eq: d1}
\coth d_{ij}=\frac{\sinh(f_j-f_i-C_{ij})}{\sinh l_{ij}}+\frac{\eta_{ij}e^{f_i+f_j}}{\sinh l_{ij}}
\end{equation*}
with
\begin{equation}\label{Eq: DCS1}
\cosh l_{ij}
=-\cosh(f_j-f_i-C_{ij})+\eta_{ij}e^{f_i+f_j},
\end{equation}
where $C\in \mathbb{R}^{E_+}$ is a constant vector satisfying
$C_{ij}+C_{jk}+C_{ki}=0$ for any $\{ijk\}\in F$ and $C_{rs}+C_{sr}=0$ for any subset $\{r,s\}\subseteq\{i,j,k\}$.
\end{description}
\item[(B)] if $\frac{\sinh d_{ij}}{\sinh d_{ji}}<0$, then the discrete conformal structure $d=d(f)$ has one of the following two forms
  \begin{description}
    \item[(i)] \begin{equation*}
\coth d_{ij}
=\frac{\alpha_ie^{2f_i}}{\sinh l_{ij}}
\sqrt{\frac{1+\alpha_je^{2f_j}}{1+\alpha_ie^{2f_i}}}
+\frac{\eta_{ij}e^{f_i+f_j}}{\sinh l_{ij}}
\end{equation*}
with
\begin{equation}\label{Eq: DCS4}
\cosh l_{ij}
=\sqrt{(1+\alpha_ie^{2f_i})(1+\alpha_je^{2f_j})}
+\eta_{ij}e^{f_i+f_j},\ \text{for}\ 1+\alpha_ie^{2f_i}>0, 1+\alpha_je^{2f_j}>0,
\end{equation}
\begin{equation}\label{Eq: new 2}
\cosh l_{ij}
=-\sqrt{(1+\alpha_ie^{2f_i})(1+\alpha_je^{2f_j})}
+\eta_{ij}e^{f_i+f_j},\ \text{for}\ 1+\alpha_ie^{2f_i}<0, 1+\alpha_je^{2f_j}<0;
\end{equation}
\item[(ii)]
\begin{equation*}
\coth d_{ij}=-\frac{\sinh(f_j-f_i-C_{ij})}{\sinh l_{ij}}+\frac{\eta_{ij}e^{f_i+f_j}}{\sinh l_{ij}}
\end{equation*}
with
\begin{equation}\label{Eq: DCS2}
\cosh l_{ij}
=\cosh(f_j-f_i-C_{ij})+\eta_{ij}e^{f_i+f_j},
\end{equation}
where $C\in \mathbb{R}^{E_+}$ is a constant vector satisfying $C_{ij}+C_{jk}+C_{ki}=0$ for any $\{ijk\}\in F$
and $C_{rs}+C_{sr}=0$ for any subset $\{r,s\}\subseteq\{i,j,k\}$.
  \end{description}
\end{description}
Furthermore, we have
\begin{description}
  \item[(a)] any one of the discrete conformal structures (\ref{Eq: DCS3}), (\ref{Eq: new 1}) or (\ref{Eq: DCS1}) can exist alone on an ideally triangulated surface with boundary, while the discrete conformal structures (\ref{Eq: DCS3}), (\ref{Eq: new 1}) and (\ref{Eq: DCS1}) can not exist simultaneously in pairs on the same ideally triangulated surface with boundary.
  \item[(b)] none of the discrete conformal structures  (\ref{Eq: DCS4}), (\ref{Eq: new 2}) and (\ref{Eq: DCS2}) can exist alone on an ideally triangulated surface with boundary, and the discrete conformal structures  (\ref{Eq: DCS4}), (\ref{Eq: new 2}) and (\ref{Eq: DCS2})
      can not exist simultaneously in pairs on the same ideally triangulated surface with boundary.
  \item[(c)] the only mixed types of discrete conformal structures on an ideally triangulated surface with boundary are  (\ref{Eq: DCS3}) (\ref{Eq: DCS4}), (\ref{Eq: new 1}) (\ref{Eq: new 2}) and (\ref{Eq: DCS1}) (\ref{Eq: DCS2}).
\end{description}
\end{theorem}

\begin{remark}
Following Thomas (\cite{Thomas}, page 53),
one can reparameterize the discrete conformal structure $d=d(f)$ so that $\alpha: B\rightarrow \{-1,0,1\}$ in (\ref{Eq: DCS3}) while keeping the induced discrete hyperbolic metric $l$ invariant.
For example, if $\alpha_i, \alpha_j>0$, then
\begin{equation*}
\begin{aligned}
\cosh l_{ij}
=&-[(1+\alpha_ie^{2f_i})(1+\alpha_je^{2f_j})]
^{1/2}+\eta_{ij}e^{f_i+f_j}\\
=&-[(1+e^{2(f_i+\frac{1}{2}\log \alpha_i)})
(1+e^{2(f_j+\frac{1}{2}\log \alpha_j)})]^{1/2}\\
&+\frac{\eta_{ij}}{(\alpha_i\alpha_j)^{1/2}}
e^{(f_i+\frac{1}{2}\log \alpha_i)+(f_j+\frac{1}{2}\log \alpha_j)}\\
=&-[(1+e^{2g_i})(1+e^{2g_j})]^{1/2}
+\widetilde{\eta}_{ij}e^{g_i+g_j},
\end{aligned}
\end{equation*}
where $g_i=f_i+\frac{1}{2}\log \alpha_i, \widetilde{\eta}_{ij}
=\frac{\eta_{ij}}{(\alpha_i\alpha_j)^{1/2}}$.
Therefore, the discrete conformal structure (\ref{Eq: DCS3}) with parameters $(\alpha_i,\alpha_j,\eta_{ij})$ is reparameterized to be  the discrete conformal structure with parameters $(1,1,\widetilde{\eta}_{ij})$.
Similar arguments also apply for the parameters in the discrete conformal structure (\ref{Eq: DCS4}).
For simplicity, we always assume $\alpha: B\rightarrow \{-1,0,1\}$ in the discrete conformal structures (\ref{Eq: DCS3}) and (\ref{Eq: DCS4}) in the following.
The conditions that $1+\alpha_ie^{2f_i}<0, 1+\alpha_je^{2f_j}<0$ in the discrete conformal structures (\ref{Eq: new 1}) and (\ref{Eq: new 2}) imply $\alpha_i<0, \alpha_j<0$.
Thus we always assume $\alpha\equiv -1$ in the discrete conformal structures (\ref{Eq: new 1}) and (\ref{Eq: new 2}) in the following.

Moreover, if $(\Sigma,\mathcal{T})$ is an ideally triangulated surface with boundary and the genus is $0$,
one can also reparameterize discrete conformal structures $d=d(f)$ in (\ref{Eq: DCS1}) and (\ref{Eq: DCS2}) so that $C\equiv0$ while keeping the induced discrete hyperbolic metric invariant.
Indeed, for these discrete conformal structures in  (\ref{Eq: DCS1}) and (\ref{Eq: DCS2}), we have
\begin{eqnarray}\label{Eq: F90}
\begin{cases}
C_{ij}+C_{jk}+C_{ki}=0,\\
C_{ij}+C_{ji}=0,\\
C_{jk}+C_{kj}=0,\\
C_{ki}+C_{ik}=0
\end{cases}
\end{eqnarray}
for any ideal face $\{ijk\}\in F$.
Fix $p_0\in B$.
By the connectedness of $(\Sigma,\mathcal{T})$, for any $p\in B$, there exists a path $\gamma_1: p_0\sim p_1\sim...\sim p_{n-1}\sim p_n=p$ connecting $p_0$ and $p$.
We define a function $g: B\rightarrow \mathbb{R}$ by setting $g(p)=g_n=-\sum_{i=1}^nC_{i-1,i}$ and $g_0=0$.
Then $g(p)$ is well-defined in the sense that it is independent of the paths connecting $p_0$ and $p$.
In fact, if there exists another path $\gamma_2: p_0\sim p^\prime_1\sim...\sim p^\prime_{m-1}\sim p^\prime_m=p$ connecting $p_0$ and $p$,
then $\widetilde{g}_m=-\sum_{k=1}^mC^\prime_{k-1,k}$ by definition.
Since the genus of $(\Sigma,\mathcal{T})$ is zero, then by (\ref{Eq: F90}), we have $g_n=\widetilde{g}_m$.
This ensures the well-definiteness of $g$.
As a result, $C_{ij}=g_{i}-g_j$ for $i$ adjacent to $j$.
We claim $g$ is unique up to a constant.
Suppose there exists another function $\overline{g}: B\rightarrow \mathbb{R}$ such that
$C_{ij}=\overline{g}_{i}-\overline{g}_j$ for $i$ adjacent to $j$.
Then $g_{i}-\overline{g}_{i}=g_j-\overline{g}_j$  for $i$ adjacent to $j$.
By the connectedness of $(\Sigma,\mathcal{T})$, we have $g-\overline{g}=c(1,...,1)$, where $c\in \mathbb{R}$ is a constant.
Therefore, (\ref{Eq: DCS1}) is equivalent to
\begin{equation*}
\cosh l_{ij}
=\cosh(f_j-f_i-g_i+g_j)+\eta_{ij}e^{f_i+f_j}
=\cosh(h_j-h_i)+\widetilde{\eta}_{ij}e^{h_i+h_j},
\end{equation*}
where $h_i=f_i+g_i$, $\widetilde{\eta}_{ij}=e^{-g_i-g_j}\eta_{ij}$.
Similar arguments also apply to the discrete conformal structure (\ref{Eq: DCS2}).
However, for a general ideally triangulated surface with boundary $(\Sigma,\mathcal{T})$, i.e., the genus may not be zero, we take $C\equiv0$ in the discrete conformal structures (\ref{Eq: DCS1}) and (\ref{Eq: DCS2}) as special cases.
\end{remark}

\begin{remark}
By Theorem \ref{Thm: DCS}, the discrete conformal structure in Definition \ref{Def: DCS} unifies and generalizes the existing most discrete conformal structures on surfaces with boundary.
More precisely, the discrete conformal structure (\ref{Eq: DCS3}) can be reduced to Guo-Luo's $(-1,-1,-1)$ type generalized circle packings \cite{GL2} if $\alpha\equiv1$, to Guo's vertex scalings \cite{Guo} if $\alpha\equiv0$ and to Xu's discrete conformal structures \cite{Xu22} if $\alpha\equiv-1$.
The discrete conformal structure (\ref{Eq: new 1}) is exactly Guo-Luo's $(-1,-1,1)$ type generalized circle packings \cite{GL2} if $\alpha\equiv-1$.
The discrete conformal structure (\ref{Eq: DCS1}) can be reduced to Guo-Luo's $(-1,-1,0)$ type generalized circle packings \cite{GL2} if $C_{ij}=0$ for any edge $\{ij\}\in E$.
Moreover, the discrete conformal structure (\ref{Eq: DCS4}) is formally the same as Glickenstein-Thomas's discrete conformal structure in the hyperbolic background geometry \cite{GT}.
However, they involve some twisted generalized hyperbolic triangles in Roger-Yang \cite{Roger-Yang}.
The discrete conformal structure (\ref{Eq: DCS4}) is closely related to the $(-1,-1,-1)$ type of twisted generalized hyperbolic triangle if $\alpha\equiv1$, to $(1,1,-1)$ type of twisted generalized hyperbolic triangle if $\alpha\equiv0$ and to the $(0,0,-1)$ type of twisted generalized hyperbolic triangle if $\alpha\equiv-1$.
The discrete conformal structure (\ref{Eq: new 2}) is closely related to the $(1,-1,-1)$ type of twisted generalized hyperbolic triangle if $\alpha\equiv-1$.
The discrete conformal structure (\ref{Eq: DCS2}) is closely related to the $(0,-1,-1)$ type of twisted generalized hyperbolic triangle if $C_{ij}=0$ for any edge $\{ij\}\in E$.
Please refer to Section \ref{section 6} for more details.
\end{remark}

\subsection{Organization of the paper}
In Section \ref{section 2}, we first give some preliminaries on the hyperboloid model of the hyperbolic space.
Then we give the motivations of Definition \ref{Def: partial edge length} and Definition \ref{Def: DCS} and prove Theorem \ref{Thm: DCS}.
In Section \ref{section 6}, we discuss the relationships between discrete conformal structures on ideally triangulated surfaces with boundary,
generalized hyperbolic triangles (including twisted generalized hyperbolic triangles) and the 3-dimensional hyperbolic geometry.
\\
\\
\textbf{Acknowledgements}\\[8pt]
This work was once presented on the workshop  ``Discrete Conformal Geometry" taking place
at Wuhan University in China from June 21st to June 24th, 2024.
During the workshop, Professor Feng Luo at Rutgers University
and Professor Tian Yang at Texas A\&M University suggested the relationships between the twisted generalized triangles and the discrete conformal structures on surfaces with boundary.
The authors thank Professor Feng Luo and Professor Tian Yang for their suggestions.
The research of the authors is supported by National Natural Science Foundation of China
under grant no. 12471057.

\section{Characterization of discrete conformal structures}\label{section 2}

In this section, we first present some basic facts in the 2-dimensional hyperbolic geometry, which are discussed in Chapter 3 and 6 of \cite{Ratcliffe}.
Then we give the motivations of Definition \ref{Def: partial edge length} and Definition \ref{Def: DCS},
the idea of which comes from Glickenstein \cite{Glickenstein} and Glickenstein-Thomas \cite{GT}.
In the end, we prove Theorem \ref{Thm: DCS}.

\subsection{Lorentzian cross product}

Let $x=(x_1,x_2,x_3),\ y=(y_1,y_2,y_3)$ be vectors in the 3-dimensional Lorentzian space $\mathbb{R}^3$.
The Lorentzian inner product $*$ of $x$ and $y$ is defined to be
\begin{equation*}
x*y=x_1y_1+x_2y_2-x_3y_3=(x_1,x_2,x_3)
\left(
   \begin{array}{ccc}
     1 & 0 & 0 \\
     0 & 1 & 0 \\
     0 & 0 & -1 \\
   \end{array}
 \right)\\
\left(
   \begin{array}{c}
     y_1  \\
     y_2  \\
     y_3  \\
   \end{array}
 \right)\\
:=xJy^T,
\end{equation*}
where $J=\mathrm{diag} \{1,1, -1\}$.
Two vectors $x,y\in \mathbb{R}^3$ are Lorentz orthogonal if and only if $x*y=0$.
The Lorentzian cross product $\otimes$ of $x$ and $y$ is defined to be
\begin{equation*}
x\otimes y: =J(x\times y),
\end{equation*}
where $\times$ is the Euclidean cross product.
Then $x\otimes y=0$ if and only if $x$ and $y$ are linearly dependent.

\begin{proposition}[\cite{Ratcliffe}, Theorem 3.2.1]\label{Prop: cross product}
If $x,y,z,w$ are vectors in $\mathbb{R}^3$, then
\begin{description}
\item[$(i)$] $x\otimes y=-y\otimes x$;
\item[$(ii)$] $(x\otimes y)\ast z=\det(x,y,z)$;
\item[$(iii)$] $x\otimes(y\otimes z)=(x*y)z-(z*x)y$;
\item[$(iv)$] $(x\otimes y)*(z\otimes w)=
\left|
   \begin{array}{cc}
     x* w & x* z \\
     y* w & y* z \\
   \end{array}
 \right|$.
\end{description}
\end{proposition}

The Lorentzian norm of $x\in \mathbb{R}^3$ is defined to be the complex number $||x||=\sqrt{x*x}$.
Here $||x||$ is either positive, zero or positive imaginary.
If $||x||$ is positive imaginary, we denote its absolute value by $|||x|||$.
Furthermore, $x\in \mathbb{R}^3$ is called space-like if $x*x>0$, light-like if $x*x=0$ and time-like if $x*x<0$.
The hyperboloid model of $\mathbb{H}^2$ is defined to be
\begin{equation*}
\mathbb{H}^2=\{x=(x_1,x_2,x_3)\in \mathbb{R}^3|\ x*x=-1\ \text{and}\ x_3>0 \},
\end{equation*}
which can be embedded in $(\mathbb{R}^3,*)$.
In fact, $(\mathbb{H}^2, *|_{\mathbb{H}^2})$ is a 2-dimensional Riemannian manifold.
A geodesic (or hyperbolic line) in $\mathbb{H}^2$ is the non-empty intersection of $\mathbb{H}^2$ with a 2-dimensional vector subspace of $\mathbb{R}^3$.

The Klein model of $\mathbb{H}^2$ is the center projection of $\mathbb{H}^2$ into the plane $z=1$ of $\mathbb{R}^3$.
In the Klein model, the time-like vectors are represented by the interior points of the unit disk, which are called as hyperbolic points.
The light-like vectors are represented by the boundary points of the unit disk, which are called as ideal points.
The space-like vectors are represented by the external points of the unit disk, which are called as hyper-ideal points.
A geodesic in $\mathbb{H}^2$ is a chord on the unit disk.
In the Klein model, a triangle $xyz$ is called a hyperbolic triangle if its three vertices $x,y,z$ are hyperbolic, and called a generalized hyperbolic triangle if at least one of its three vertices is ideal or hyper-ideal.
If the vertices of a generalized hyperbolic triangle are all hyper-ideal, we call this triangle as a hyper-ideal hyperbolic triangle for simplicity.

As a direct application of the identity $(iv)$ in Proposition \ref{Prop: cross product}, we have the following result.

\begin{proposition}\label{Prop: Right angle}
Let $xyz$ be a (hyperbolic or generalized hyperbolic) triangle with $x\in\mathbb{H}^2$ and a right angle at $x$. Then
\begin{equation*}
-(z*y)=(z*x)(x*y).
\end{equation*}
\end{proposition}
\proof
Direct calculations give
\begin{align*}
0=(x\otimes y)*(x\otimes z)
=&\left|
   \begin{array}{cc}
     x* z & x* x \\
     y* z & y* x \\
   \end{array}
 \right|\\
=&(z*x)(x*y)-(z*y)(x* x)\\
=&(z*x)(x*y)+(z*y).
\end{align*}
\qed

For any two vectors $x,y\in \mathbb{H}^2$, there is a unique geodesic through $x$ and $y$.
The arc length of the geodesic between $x$ and $y$ is defined to be the hyperbolic distance between $x$ and $y$, denoted by $d_{\mathbb{H}}(x,y)$.
We use $\mathrm{Span}(x)$ and $\mathrm{Span}(x,y)$ to
denote the vector subspace spanned $x$ and $x,y$, respectively.
If $x\in \mathbb{R}^3$ is a space-like vector,
then the Lorentzian complement $x^{\bot}:=\{w\in \mathbb{R}^3| x*w=0\}$ of $\mathrm{Span}(x)$ is a 2-dimension vector subspace of $\mathbb{R}^3$ and $x^{\bot}\cap \mathbb{H}^2$ is non-empty.
This implies that $x^{\bot}\cap \mathbb{H}^2$ is a geodesic in $\mathbb{H}^2$.
And the hyperbolic distance between $x$ and $y\in \mathbb{H}^2$ is defined to be $d_{\mathbb{H}}(y,x^\bot)
=\inf\{d_{\mathbb{H}}(y,z)|z\in x^\bot\cap \mathbb{H}^2\}$.
Let $x,y\in \mathbb{R}^3$ be space-like vectors.
If $\mathrm{Span}(x,y)\cap \mathbb{H}^2$ is non-empty, then the hyperbolic distance between $x$ and $y$ is defined by $d_{\mathbb{H}}(x^\bot,y^\bot)
=\inf\{d_{\mathbb{H}}(z,w)|z\in x^\bot\cap \mathbb{H}^2, w\in y^\bot\cap \mathbb{H}^2\}$.
If $\mathrm{Span}(x,y)\cap \mathbb{H}^2$ is empty, then $x^\bot$ and $y^\bot$ intersect in an angle $\angle(x^\bot,y^\bot)$ in $\mathbb{H}^2$.

\begin{proposition}[\cite{Ratcliffe}, Chapter 3]\label{Prop: Ratcliffe 1}
Suppose $x,y\in\mathbb{R}^3$ are time-like vectors and $z,w\in \mathbb{R}^3$ are space-like vectors. Then
\begin{description}
\item[$(i)$] $x*y=-|||x|||\cdot|||y|||\cosh d_{\mathbb{H}}(x,y)$;
\item[$(ii)$] $|x*z|=|||x|||\cdot||z||\sinh d_{\mathbb{H}}(x,z^\bot)$, and $x*z<0$ if and only if $x$ and $z$ are on opposite sides of the hyperplane $z^\bot$;
\item[$(iii)$] if $\mathrm{Span}(z,w)\cap \mathbb{H}^2$ is non-empty, then $|z*w|=||z||\cdot||w||\cosh d_{\mathbb{H}}(z^\bot,w^\bot)$, and $z*w<0$ if and only if $z$ and $w$ are oppositely oriented tangent vectors of $Q$, where $Q$ is the hyperbolic line Lorentz orthogonal to $z^\bot$ and $w^\bot$;
\item[$(iv)$] if $\mathrm{Span}(z,w)\cap \mathbb{H}^2$ is empty, then $z*w=||z||\cdot||w||\cos \angle(z^\bot,w^\bot)$.
\end{description}
\end{proposition}

\subsection{Compatible condition}
As mentioned in \cite{Thomas}, not every assignment of $d\in \mathbb{R}^{E_+}$ with $l_{ij}=d_{ij}+d_{ji}>0$ induces a good geometric structure on the Poincar\'{e} dual, because the three perpendiculars of a right-angled hyperbolic hexagons may fail to meet at a common point in the Klein model.
To ensure that there exists a good geometric structure on the Poincar\'{e} dual,
we require that the three perpendiculars $P_{ij}, P_{jk}, P_{ki}$ of a right-angled hyperbolic  hexagon $\{ijk\}\in F$ intersect at a common point.
The following lemma motivates Definition \ref{Def: partial edge length}.

\begin{lemma}\label{Lem: compatible condition}
Suppose $\{ijk\}$ is a right-angled hyperbolic  hexagon with $d\in \mathbb{R}^{E_+}$ satisfying $l_{rs}=d_{rs}+d_{sr}$.
Then the three perpendiculars $P_{ij}, P_{jk}, P_{ki}$ intersect at a common point if and only if the following compatible formula
\begin{equation}\label{Eq: compatible condition 3}
(c_{ij}\ast p_i)(c_{jk}\ast p_j)(c_{ki}\ast p_k)=
(c_{ij}\ast p_j)(c_{jk}\ast p_k)(c_{ki}\ast p_i)
\end{equation}
holds, which is equivalent to the formula
$$\sinh d_{ij} \sinh d_{jk} \sinh d_{ki}=\sinh d_{ji} \sinh d_{kj} \sinh d_{ik}$$
in (\ref{Eq: compatible condition}).
\end{lemma}
\proof
The trick is taken from Glickenstein-Thomas \cite{GT}.
Without loss of generality, we assume $||v_i||=||v_j||=||v_k||=1$.
If $c\in \mathbb{H}^2$ is any other point in $P_{ij}$ and linearly independent of edge center $c_{ij}$, then $P_{ij}=(c\otimes c_{ij})^\bot\cap\mathbb{H}^2$.
Note that $E_{ij}=\mathrm{Span}(v_i,v_j)\cap \mathbb{H}^2$
can be written as $(v_i\otimes v_j)^\bot\cap \mathbb{H}^2$.
Since $P_{ij}$ and $E_{ij}$ are Lorentz orthogonal,
then
\begin{equation*}
(c\otimes c_{ij})\ast(v_i\otimes v_j)=0.
\end{equation*}
Using the identity $(iv)$ in Proposition \ref{Prop: cross product}, we have
\begin{equation*}
c\ast[(c_{ij}\ast v_i)v_j-(c_{ij}\ast v_j)v_i]=0.
\end{equation*}
Therefore, the three perpendiculars $P_{ij}, P_{jk}, P_{ki}$ intersect in a common point if and only if there exists a non-trivial solution $c$ satisfying the following system
\begin{equation*}
\begin{aligned}
&c\ast[(c_{ij}\ast v_i)v_j-(c_{ij}\ast v_j)v_i]=0,\\
&c\ast[(c_{jk}\ast v_j)v_k-(c_{jk}\ast v_k)v_j]=0,\\
&c\ast[(c_{ki}\ast v_k)v_i-(c_{ki}\ast v_i)v_k]=0.
\end{aligned}
\end{equation*}
This system can be written in the following matrix form
\begin{equation}\label{Eq: equation for cijk}
\left[ {
\begin{matrix}
  [(c_{ij}\ast v_i)v_j-(c_{ij}\ast v_j)v_i]^T  \\
  [(c_{jk}\ast v_j)v_k-(c_{jk}\ast v_k)v_j]^T  \\
  [(c_{ki}\ast v_k)v_i-(c_{ki}\ast v_i)v_k]^T   \\
\end{matrix} }\right]
\cdot J\cdot c:=M\cdot J\cdot c=0.
\end{equation}
The system (\ref{Eq: equation for cijk}) has a non-trivial solution
if and only if the determinant of the matrix $M$ is zero because of the non-singularity of $J$.
A routine calculation gives rise to
\begin{equation*}
\det(M)=[(c_{ij}\ast v_i)(c_{jk}\ast v_j)(c_{ki}\ast v_k)-
(c_{ij}\ast v_j)(c_{jk}\ast v_k)(c_{ki}\ast v_i)]\det(v_i, v_j, v_k).
\end{equation*}
This implies the formula (\ref{Eq: compatible condition 3}) holds by $\det(v_i, v_j, v_k)\neq0$.
Since the edge center $c_{ij}$ is time-like of length $-1$ and $v_i$ is space-like of length $1$, then $c_{ij}*v_i=-\sinh d_{ij}$ by Proposition \ref{Prop: Ratcliffe 1}.
This implies (\ref{Eq: compatible condition 3}) is equivalent to (\ref{Eq: compatible condition}).
\qed

\subsection{Conformal variations of a right-angled hyperbolic  hexagon}\label{subsection 2}

\begin{figure}[!ht]
\centering
\includegraphics[scale=1]{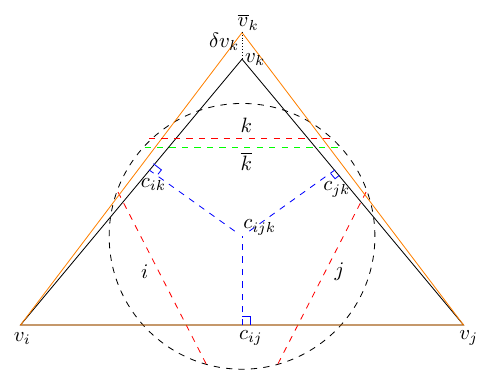}
\caption{A hyper-ideal hyperbolic triangle $v_iv_jv_k$ under $f_k-$conformal variation in the Klein model. }
\label{Figure_5}
\end{figure}

For a right-angled hyperbolic hexagon $\{ijk\}\in F$, there exists a unique hyper-ideal hyperbolic triangle $\{ijk\}$ corresponding to it with the vertices corresponding to the boundary arcs.
We consider the effect of perturbing $f_k$ conformally while keeping $f_i$ and $f_j$ fixed, i.e., $\delta f_i=\delta f_j=0$. This is called an $f_k-$conformal variation by Glickenstein-Thomas \cite{GT}.
In \cite{GT}, Glickenstein-Thomas added a condition to a hyperbolic triangle in a closed surface, requiring that the vertex $v_k$ and face center $c_{ijk}$ of a hyperbolic triangle should be in a line including the vertex $\overline{v}_k$ of the new hyperbolic triangle under the $f_k-$conformal variation.
We change the condition ``in a line'' in Glickenstein-Thomas \cite{GT} to the condition ``in a 2-dimensional vector subspace of $\mathbb{R}^3$'' for a hyper-ideal hyperbolic triangle in an ideally triangulated surface with boundary.
In the Klein model, this is almost equivalent to the condition introduced by Glickenstein-Thomas, except that the vertices $v_k$ and $\overline{v}_k$ are spacelike.
Please refer to Figure \ref{Figure_5}.
Note that if the vertex $v_i$ of a hyper-ideal hyperbolic triangle $v_iv_jv_k$ is fixed, then the position of boundary component $i=v_i^\bot\cap\mathbb{H}^2$ will not move, while the arc length may change.

The following lemma motivates Definition \ref{Def: DCS}.
\begin{lemma}\label{Lem: conformal variation}
Let $\{ijk\}$ be a right-angled hyperbolic hexagon with partial edge lengths $d\in \mathbb{R}^{E_+}$.
Under the $f_k-$conformal variation,
the vectors $\overline{v}_k, v_k$ and $c_{ijk}$ lie in a 2-dimensional vector subspace of the Lorentzian space  $\mathbb{R}^3$ if and only if there exists a function $G: \mathbb{R}\rightarrow \mathbb{R}$ satisfying
$$\frac{\partial l_{sk}}{\partial f_k}=G(f_k)\coth d_{ks}$$
for $s=i,j$.
\end{lemma}
\proof
For simplicity, set $\{1,2,3\}=\{i,j,k\}$.
Without loss of generality, we assume $||v_1||=||v_2||=||v_3||=1$ and $c_{123}*c_{123}=\pm1$. The case that $c_{123}$ is lightlike follows from continuity by moving the face center $c_{123}$ to the light cone.

Using the bilinearity of Lorentzian inner product gives
\begin{equation}\label{Eq: F1}
v_s\ast\delta v_3=\delta(v_s\ast v_3)=\delta(-\cosh l_{s3})=-\sinh l_{s3}\cdot\frac{\partial l_{s3}}{\partial f_3}\cdot\delta f_3
\end{equation}
for $s=1,2$.
The vectors $\overline{v}_3, v_3$ and $c_{123}$ lie in a 2-dimensional vector subspace of $\mathbb{R}^3$ if and only if there exist two constants $\lambda, \mu\in \mathbb{R}$ such that
\begin{equation}\label{Eq: F84}
\overline{v}_3=\lambda v_3+\mu c_{123}.
\end{equation}
Since
\begin{equation*}
0=\frac{1}{2}\delta(v_3\ast v_3)=v_3\ast\delta v_3=v_3\ast(\overline{v}_3-v_3)=v_3\ast \overline{v}_3-v_3\ast v_3=v_3\ast \overline{v}_3-1,
\end{equation*}
then
\begin{equation*}
1=v_3\ast \overline{v}_3=v_3\ast(\lambda v_3+\mu c_{123})=\lambda+\mu(v_3\ast c_{123}),
\end{equation*}
which implies
\begin{equation}\label{Eq: lambda}
\lambda=1-\mu(v_3\ast c_{123}).
\end{equation}
Set $\mu=-(\delta f_3)t$.
Then by (\ref{Eq: F84}) and (\ref{Eq: lambda}), we have
\begin{equation*}
\overline{v}_3
=\lambda v_3-(\delta f_3)t c_{123}
=[1+(v_3\ast c_{123})(\delta f_3)t]v_3
-(\delta f_3)t c_{123},
\end{equation*}
which implies
\begin{equation*}
\delta v_3=\overline{v}_3-v_3
=(v_3\ast c_{123})(\delta f_3)t v_3-(\delta f_3)t c_{123}.
\end{equation*}
Therefore,
\begin{equation}\label{Eq: F2}
\begin{aligned}
v_1\ast \delta v_3
&=(v_1\ast v_3)(v_3\ast c_{123})(\delta f_3)t-(v_1\ast c_{123})(\delta f_3)t\\
&=(v_3\ast c_{123})(\delta f_3)t\left[(v_1\ast v_3)-\frac{v_1\ast c_{123}}{v_3\ast c_{123}}\right].
\end{aligned}
\end{equation}
By Proposition \ref{Prop: Right angle},
in the generalized right-angled hyperbolic
triangle $v_3c_{13}c_{123}$,
we have
\begin{equation}\label{Eq: v3a}
v_3\ast c_{123}
=-(v_3\ast c_{13})(c_{13}\ast c_{123})
=\sinh d_{31}\cdot(c_{13}\ast c_{123}).
\end{equation}
Similarly, in the generalized right-angled hyperbolic  triangle $v_1c_{13}c_{123}$, we have
\begin{equation}\label{Eq: v3b}
v_1\ast c_{123}
=-(v_1\ast c_{13})(c_{13}\ast c_{123})
=\sinh d_{13}\cdot(c_{13}\ast c_{123}).
\end{equation}
Thus by (\ref{Eq: v3a}) and  (\ref{Eq: v3b}), we have
\begin{equation}\label{Eq: F3}
\frac{v_1\ast c_{123}}{v_3\ast c_{123}}
=\frac{\sinh d_{13}}{\sinh d_{31}}.
\end{equation}
Substituting $v_1\ast v_3=-\cosh l_{13}$ and (\ref{Eq: F3}) into (\ref{Eq: F2}) gives
\begin{equation}\label{Eq: F4}
\begin{aligned}
v_1\ast \delta v_3
=&(v_3\ast c_{123})(\delta f_3)t \left[-\cosh l_{13}-\frac{\sinh d_{13}}{\sinh d_{31}} \right]\\
=&-(v_3\ast c_{123})(\delta f_3)t\cdot\frac{1}{\sinh d_{31}}\bigg(\cosh l_{13}\sinh d_{31}+\sinh(l_{13}-d_{31})\bigg)\\
=&-(v_3\ast c_{123})(\delta f_3)t\sinh l_{13}\coth d_{31}.
\end{aligned}
\end{equation}
Similarly,
\begin{equation}\label{Eq: F5}
v_2\ast \delta v_3
=-(v_3\ast c_{123})(\delta f_3)t\sinh l_{23}\coth d_{32}.
\end{equation}
Combining (\ref{Eq: F1}), (\ref{Eq: F4}) and (\ref{Eq: F5}), we have
\begin{equation*}
\frac{\partial l_{s3}}{\partial f_3}
=(v_3\ast c_{123})t\coth d_{3s}
:=\widetilde{\lambda}\coth d_{3s}
\end{equation*}
for $s=1,2$,
which implies $\widetilde{\lambda}=\tanh d_{3s}\cdot\frac{\partial l_{s3}}{\partial f_3}$.
Note that $l_{s3}=d_{s3}+d_{3s}$ depends only on the parameters $f_s$ and $f_3$ by (\ref{Eq: variation 2}), then $\widetilde{\lambda}$ depends only on the parameters $f_s$ and $f_3$.
Therefore, there exist functions $G_1: \mathbb{R}\rightarrow \mathbb{R}$ and $G_2: \mathbb{R}\rightarrow \mathbb{R}$ such that
$\widetilde{\lambda}=G_1(f_1,f_3)=G_2(f_2,f_3)$, which is impossible unless $\widetilde{\lambda}=G(f_3)$.
This completes the proof.
\qed

\begin{remark}
By the proof of Lemma \ref{Lem: conformal variation}, the function $G$ depends only on $f_k$.
By a change of variables, we can set $G(f_k)\equiv1$ for simplicity, which gives Definition \ref{Def: DCS}.
\end{remark}

\subsection{Proof of Theorem \ref{Thm: DCS}}
\label{subsection 1}
By the formula (\ref{Eq: variation 1}), we have
\begin{equation}\label{Eq: F28}
\frac{\partial}{\partial f_i}\cosh l_{ij}
=\sinh l_{ij}\frac{\partial l_{ij}}{\partial f_i}
=\sinh l_{ij}\coth d_{ij}
=\cosh l_{ij}+\frac{\sinh d_{ji}}{\sinh d_{ij}}.
\end{equation}
Similarly,
\begin{equation}\label{Eq: F29}
\frac{\partial}{\partial f_j}\cosh l_{ij}
=\cosh l_{ij}+\frac{\sinh d_{ij}}{\sinh d_{ji}}.
\end{equation}
Note that
\begin{equation*}
\frac{\partial}{\partial f_j}\coth d_{ij}
=\frac{\partial^2 l_{ij}}{\partial f_i\partial f_j}
=\frac{\partial}{\partial f_i}\coth d_{ji},
\end{equation*}
which is equivalent to
\begin{equation}\label{Eq: F20}
\frac{1}{\sinh^2 d_{ij}}\frac{\partial d_{ij}}{\partial f_j}=\frac{1}{\sinh^2 d_{ji}}\frac{\partial d_{ji}}{\partial f_i}.
\end{equation}
Set
\begin{equation}\label{Definition of H}
H=\log\frac{\sinh^2 d_{ij}}{\sinh^2 d_{ji}}.
\end{equation}
Then the formulas (\ref{Eq: compatible condition}) and (\ref{Eq: variation 2}) give
\begin{equation}\label{Eq: F21}
\frac{\partial^2 H}{\partial f_i\partial f_j}=0.
\end{equation}
Direct calculations give
\begin{equation*}
\begin{aligned}
\frac{\partial H}{\partial f_i}
=&2\frac{\sinh d_{ji}}{\sinh d_{ij}}\frac{\partial}{\partial f_i}(\frac{\sinh d_{ij}}{\sinh d_{ji}})\\
=&2\frac{\sinh d_{ji}}{\sinh d_{ij}}\frac{1}{\sinh^2 d_{ji}}(\cosh d_{ij}\frac{\partial d_{ij}}{\partial f_i}\sinh d_{ji}-\cosh d_{ji}\frac{\partial d_{ji}}{\partial f_i}\sinh d_{ij})\\
=&2(\coth d_{ij}\frac{\partial d_{ij}}{\partial f_i}-\coth d_{ji}\frac{\partial d_{ji}}{\partial f_i})\\
=&2\left[\coth d_{ij}\frac{\partial d_{ij}}{\partial f_i}-\coth d_{ji}\bigg(\frac{\sinh^2 d_{ji}}{\sinh^2 d_{ij}}\bigg)\frac{\partial d_{ij}}{\partial f_j}\right],
\end{aligned}
\end{equation*}
where (\ref{Eq: F20}) is used in the last line.
Similarly,
\begin{equation*}
\begin{aligned}
\frac{\partial H}{\partial f_j}
=&2(\coth d_{ij}\frac{\partial d_{ij}}{\partial f_j}-\coth d_{ji}\frac{\partial d_{ji}}{\partial f_j})\\
=&2\left[\coth d_{ij}\bigg(\frac{\sinh^2 d_{ij}}{\sinh^2 d_{ji}}\bigg)\frac{\partial d_{ji}}{\partial f_i}-\coth d_{ji}\frac{\partial d_{ji}}{\partial f_j}\right],
\end{aligned}
\end{equation*}
where (\ref{Eq: F20}) is used in the last line.
Thus
\begin{equation*}
\begin{aligned}
\bigg(\frac{\sinh^2 d_{ij}}{\sinh^2 d_{ji}}\frac{\partial}{\partial f_i}+\frac{\partial}{\partial f_j}\bigg)H
=&2\left[\coth d_{ij}\bigg(\frac{\sinh^2 d_{ij}}{\sinh^2 d_{ji}}\bigg)\frac{\partial l_{ij}}{\partial f_i}-\coth d_{ji}\frac{\partial l_{ij}}{\partial f_j}\right]\\
=&2\left[\coth^2 d_{ij}\bigg(\frac{\sinh^2 d_{ij}}{\sinh^2 d_{ji}}\bigg)-\coth^2 d_{ji}\right]\\
=&2\bigg(\frac{\sinh^2 d_{ij}}{\sinh^2 d_{ji}}-1\bigg),
\end{aligned}
\end{equation*}
which is equivalent to
\begin{equation}\label{Eq: F22}
\bigg(e^H\frac{\partial}{\partial f_i}+\frac{\partial}{\partial f_j}\bigg)H=2(e^H-1),
\end{equation}
and
\begin{equation}\label{Eq: F23}
\bigg(\frac{\partial}{\partial f_i}+e^{-H}\frac{\partial}{\partial f_j}\bigg)H=2(1-e^{-H}).
\end{equation}
Using (\ref{Eq: F21}) and differentiating (\ref{Eq: F22}) with respect to $f_i$ and (\ref{Eq: F23}) with respect to $f_j$ can give
\begin{equation*}
\frac{\partial^2 H}{\partial f^2_i}+\bigg(\frac{\partial H}{\partial f_i}\bigg)^2=2\frac{\partial H}{\partial f_i},
\end{equation*}
and
\begin{equation*}
\frac{\partial^2 H}{\partial f^2_j}-\bigg(\frac{\partial H}{\partial f_j}\bigg)^2=2\frac{\partial H}{\partial f_j}.
\end{equation*}
One can easily solve this ODE to obtain that
\begin{equation*}
\frac{\partial H}{\partial f_i}\equiv2\ \quad \text{or}\ \quad
\frac{\partial H}{\partial f_i}=2\frac{a_{ij}e^{2f_i}}{1+a_{ij}e^{2f_i}}
\end{equation*}
for some constant $a_{ij}\in \mathbb{R}$ by (\ref{Eq: F21}), and
\begin{equation*}
\frac{\partial H}{\partial f_j}\equiv -2\ \quad \text{or}\ \quad
\frac{\partial H}{\partial f_j}=-2\frac{a_{ji}e^{2f_j}}{1+a_{ji}e^{2f_j}}
\end{equation*}
for some constant $a_{ji}\in \mathbb{R}$ by (\ref{Eq: F21}).
Note that the formula (\ref{Eq: F22}) can be rewritten as
\begin{equation}\label{Eq: F13}
(\frac{\partial H}{\partial f_i}-2)e^H
=-(\frac{\partial H}{\partial f_j}+2).
\end{equation}
\begin{description}
\item[(I)] If $\frac{\partial H}{\partial f_i}\equiv 2$, then $\frac{\partial H}{\partial f_j}\equiv -2$.
This implies
\begin{equation*}
H=2f_i-2f_j+c_{ij},
\end{equation*}
where $c_{ij}$ is a constant. By the definition of $H$ in (\ref{Definition of H}), we have $c_{ij}+c_{ji}=0$.
Then
\begin{equation}\label{Eq: F85}
\frac{\sinh^2 d_{ij}}{\sinh^2 d_{ji}}
=e^H=c_2e^{2f_i-2f_j},
\end{equation}
where $c_2=e^{c_{ij}}>0$.
\begin{description}
\item[(i)] If
$\frac{\sinh d_{ij}}{\sinh d_{ji}}=c_3e^{f_i-f_j}$ and
$\frac{\sinh d_{ji}}{\sinh d_{ij}}=c_4e^{f_j-f_i}$,
where $c_3=\sqrt{c_2}=e^{\frac{1}{2}c_{ij}}>0$ and $c_4=c_3^{-1}=e^{-\frac{1}{2}c_{ij}}=e^{\frac{1}{2}c_{ji}}>0$.
Then by (\ref{Eq: F28}), we have
\begin{equation*}
\frac{\partial}{\partial f_i}\cosh l_{ij}
=\cosh l_{ij}+c_4e^{f_j-f_i},
\end{equation*}
which implies
$\cosh l_{ij}=-\frac{1}{2}c_4 e^{f_j-f_i}+c_5(f_j)e^{f_i}$.
Similarly, by (\ref{Eq: F29}), we have
\begin{equation*}
\frac{\partial}{\partial f_j}\cosh l_{ij}
=\cosh l_{ij}+c_3e^{f_i-f_j},
\end{equation*}
which implies
$\cosh l_{ij}=-\frac{1}{2}c_3 e^{f_i-f_j}+c_6(f_i)e^{f_j}$.
Hence,
\begin{equation*}
-\frac{1}{2}c_4 e^{f_j-f_i}+c_5(f_j)e^{f_i}
=-\frac{1}{2}c_3 e^{f_i-f_j}+c_6(f_i)e^{f_j},
\end{equation*}
which implies
\begin{equation}\label{Eq: fi fj}
c_5(f_j)e^{-f_j}+\frac{1}{2}c_3 e^{-2f_j}
=c_6(f_i)e^{-f_i}+\frac{1}{2}c_4 e^{-2f_i}.
\end{equation}
Note that the left hand side of the equation (\ref{Eq: fi fj}) depends only on $f_j$, while
the right hand side of the equation (\ref{Eq: fi fj}) depends only on $f_i$, this implies that
the term in the equation (\ref{Eq: fi fj}) is a constant independent of $f_i$ and $f_j$, denoted by $\eta_{ij}$.
Then
$c_5(f_j)=\eta_{ij}e^{f_j}-\frac{1}{2}c_3 e^{-f_j}$,
$c_6(f_i)=\eta_{ij}e^{f_i}-\frac{1}{2}c_4 e^{-f_i}$.
Therefore,
\begin{equation*}
\cosh l_{ij}=-\cosh (f_j-f_i-C_{ij})+\eta_{ij}e^{f_i+f_j},
\end{equation*}
where $C_{ij}=\log c_3=\frac{1}{2}c_{ij}$.
Moreover, by (\ref{Eq: compatible condition}), we have
\begin{equation}\label{Eq: F6}
\log\frac{\sinh^2 d_{ij}}{\sinh^2 d_{ji}}
+\log\frac{\sinh^2 d_{jk}}{\sinh^2 d_{kj}}
+\log\frac{\sinh^2 d_{ki}}{\sinh^2 d_{ik}}
=0.
\end{equation}
Differentiating (\ref{Eq: F6}) with respect to $f_i$ gives $\frac{\partial }{\partial f_i}(\log\frac{\sinh^2 d_{ki}}{\sinh^2 d_{ik}})=-2$.
Combining with (\ref{Eq: F13}) gives $\frac{\partial }{\partial f_k}(\log\frac{\sinh^2 d_{ki}}{\sinh^2 d_{ik}})=2$.
Similarly,
$\frac{\partial }{\partial f_j}(\log\frac{\sinh^2 d_{jk}}{\sinh^2 d_{kj}})=2$ and $\frac{\partial }{\partial f_k}(\log\frac{\sinh^2 d_{jk}}{\sinh^2 d_{kj}})=-2$.
Then $c_{ij}+c_{jk}+c_{ki}=0$ by (\ref{Eq: F85}).
Hence, $C_{ij}+C_{jk}+C_{ki}=0$ and $C_{rs}+C_{sr}=0$ for $\{r,s\}\subseteq\{i,j,k\}$.
This gives the discrete conformal structure (\ref{Eq: DCS1}).
\item[(ii)]
If $\frac{\sinh d_{ij}}{\sinh d_{ji}}=-c_3e^{f_i-f_j}$, then $\frac{\sinh d_{ji}}{\sinh d_{ij}}=-c_4e^{f_j-f_i}$.
Thus by (\ref{Eq: F28}), we have
\begin{equation*}
\frac{\partial}{\partial f_i}\cosh l_{ij}
=\cosh l_{ij}-c_4e^{f_j-f_i},
\end{equation*}
which implies
$\cosh l_{ij}=\frac{1}{2}c_4 e^{f_j-f_i}+c_5(f_j)e^{f_i}$.
Similarly, by (\ref{Eq: F29}), we have
\begin{equation*}
\frac{\partial}{\partial f_j}\cosh l_{ij}
=\cosh l_{ij}-c_3e^{f_i-f_j},
\end{equation*}
which implies
$\cosh l_{ij}=\frac{1}{2}c_3 e^{f_i-f_j}+c_{6}(f_i)e^{f_j}$.
Hence,
\begin{equation*}
\frac{1}{2}c_4 e^{f_j-f_i}+c_5(f_j)e^{f_i}
=\frac{1}{2}c_3 e^{f_i-f_j}+c_6(f_i)e^{f_j},
\end{equation*}
which implies
\begin{equation}\label{Eq: fi fj b}
c_5(f_j)e^{-f_j}-\frac{1}{2}c_3 e^{-2f_j}
=c_6(f_i)e^{-f_i}-\frac{1}{2}c_4 e^{-2f_i}.
\end{equation}
Note that the left hand side and right hand side of (\ref{Eq: fi fj b}) depend only on $f_j$ and $f_i$ respectively,
the term in (\ref{Eq: fi fj b}) is a constant independent of $f_j$ and $f_i$, denoted by $\eta_{ij}$.
Then
$c_5(f_j)=\eta_{ij}e^{f_j}+\frac{1}{2}c_3 e^{-f_j}$,
$c_6(f_i)=\eta_{ij}e^{f_i}+\frac{1}{2}c_4 e^{-f_i}$.
Therefore
\begin{equation*}
\cosh l_{ij}=\cosh (f_j-f_i-C_{ij})+\eta_{ij}e^{f_i+f_j},
\end{equation*}
where $C_{ij}=\log c_3=\frac{1}{2}c_{ij}$.
Similarly, we have
$C_{ij}+C_{jk}+C_{ki}=0$ and $C_{rs}+C_{sr}=0$ for $\{r,s\}\subseteq\{i,j,k\}$.
This gives the discrete conformal structure (\ref{Eq: DCS2}).
\end{description}

\item[(II)] If $\frac{\partial H}{\partial f_i}\not\equiv 2$, by (\ref{Eq: F13}), we have
\begin{equation*}
\frac{\sinh^2 d_{ij}}{\sinh^2 d_{ji}}
=e^H=\frac{1}{\frac{\partial H}{\partial f_i}-2}(-\frac{\partial H}{\partial f_j}-2)=\frac{1+a_{ij}e^{2f_i}}{1+a_{ji}e^{2f_j}},
\end{equation*}
which implies
\begin{equation*}
\frac{\sinh d_{ij}}{\sinh d_{ji}}
=\pm\sqrt{\frac{1+a_{ij}e^{2f_i}}{1+a_{ji}e^{2f_j}}}.
\end{equation*}
\begin{description}
\item[(i)]
If $\frac{\sinh d_{ij}}{\sinh d_{ji}}
=\sqrt{\frac{1+a_{ij}e^{2f_i}}{1+a_{ji}e^{2f_j}}}$,
then by (\ref{Eq: F28}) and (\ref{Eq: F29}), one can obtain the following system
\begin{equation*}
\begin{aligned}
&\frac{\partial}{\partial f_i}\cosh l_{ij}
=\cosh l_{ij}
+\sqrt{\frac{1+a_{ji}e^{2f_j}}{1+a_{ij}e^{2f_i}}},\\
&\frac{\partial}{\partial f_j}\cosh l_{ij}
=\cosh l_{ij}
+\sqrt{\frac{1+a_{ij}e^{2f_i}}{1+a_{ji}e^{2f_j}}}.
\end{aligned}
\end{equation*}
This implies
\begin{gather*}
\cosh l_{ij}
=-\sqrt{(1+a_{ij}e^{2f_i})(1+a_{ji}e^{2f_j})}
+\eta_{ij}e^{f_i+f_j},\ \text{if}\ 1+a_{ij}e^{2f_i}>0,\ 1+a_{ji}e^{2f_j}>0, \\
\cosh l_{ij}
=\sqrt{(1+a_{ij}e^{2f_i})(1+a_{ji}e^{2f_j})}
+\eta_{ij}e^{f_i+f_j},\ \text{if}\ 1+a_{ij}e^{2f_i}<0,\ 1+a_{ji}e^{2f_j}<0, \\
\end{gather*}
where $\eta_{ij}\in \mathbb{R}$ is a constant.
By (\ref{Eq: compatible condition}) and (\ref{Eq: F21}), we have
\begin{equation}\label{Eq: F31}
\log\left|\frac{\sinh d_{ij}}{\sinh d_{ji}}\right|
+\log\left|\frac{\sinh d_{ki}}{\sinh d_{ik}}\right|
=\log\left|\frac{\sinh d_{kj}}{\sinh d_{jk}}\right|.
\end{equation}
Note that $\log|\frac{\sinh d_{kj}}{\sinh d_{jk}}|$ is independent of $f_i$ by (\ref{Eq: variation 2}), so differentiating (\ref{Eq: F31}) with respect to $f_i$ gives $a_{ij}=a_{ik}$.
Set
$$\alpha_i:=a_{ij}=a_{ik}.$$
Then
\begin{gather*}
\cosh l_{ij}
=-\sqrt{(1+\alpha_ie^{2f_i})(1+\alpha_je^{2f_j})}
+\eta_{ij}e^{f_i+f_j},\ \text{if}\ 1+\alpha_ie^{2f_i}>0,\ 1+\alpha_je^{2f_j}>0, \\
\cosh l_{ij}
=\sqrt{(1+\alpha_ie^{2f_i})(1+\alpha_je^{2f_j})}
+\eta_{ij}e^{f_i+f_j},\ \text{if}\ 1+\alpha_ie^{2f_i}<0,\ 1+\alpha_je^{2f_j}<0,
\end{gather*}
with
\begin{equation*}
\coth d_{ij}
=\frac{1}{\sinh l_{ij}}\frac{\partial}{\partial f_i}\cosh l_{ij}
=-\frac{\alpha_ie^{2f_i}}{\sinh l_{ij}}
\sqrt{\frac{1+\alpha_je^{2f_j}}{1+\alpha_ie^{2f_i}}}
+\frac{\eta_{ij}e^{f_i+f_j}}{\sinh l_{ij}}.
\end{equation*}
This gives the discrete conformal structures (\ref{Eq: DCS3}) and (\ref{Eq: new 1}).

\item[(ii)]
If $\frac{\sinh d_{ij}}{\sinh d_{ji}}
=-\sqrt{\frac{1+a_{ij}e^{2f_i}}{1+a_{ji}e^{2f_j}}}$,
then one can obtain the following system
\begin{equation*}
\begin{aligned}
&\frac{\partial}{\partial f_i}\cosh l_{ij}
=\cosh l_{ij}
-\sqrt{\frac{1+a_{ji}e^{2f_j}}{1+a_{ij}e^{2f_i}}},\\
&\frac{\partial}{\partial f_j}\cosh l_{ij}
=\cosh l_{ij}
-\sqrt{\frac{1+a_{ij}e^{2f_i}}{1+a_{ji}e^{2f_j}}}.
\end{aligned}
\end{equation*}
This implies
\begin{gather*}
\cosh l_{ij}
=\sqrt{(1+a_{ij}e^{2f_i})(1+a_{ji}e^{2f_j})}
+\eta_{ij}e^{f_i+f_j},\ \text{if}\ 1+a_{ij}e^{2f_i}>0,\ 1+a_{ji}e^{2f_j}>0, \\
\cosh l_{ij}
=-\sqrt{(1+a_{ij}e^{2f_i})(1+a_{ji}e^{2f_j})}
+\eta_{ij}e^{f_i+f_j},\ \text{if}\ 1+a_{ij}e^{2f_i}<0,\ 1+a_{ji}e^{2f_j}<0, \\
\end{gather*}
where $\eta_{ij}\in \mathbb{R}$ is a constant.
Similarly, one can also prove $a_{ij}=a_{ik}$.
Set $$\alpha_i:=a_{ij}=a_{ik}.$$
Then
\begin{gather*}
\cosh l_{ij}
=\sqrt{(1+\alpha_ie^{2f_i})(1+\alpha_je^{2f_j})}
+\eta_{ij}e^{f_i+f_j},\ \text{if}\ 1+\alpha_ie^{2f_i}>0,\ 1+\alpha_je^{2f_j}>0, \\
\cosh l_{ij}
=-\sqrt{(1+\alpha_ie^{2f_i})(1+\alpha_je^{2f_j})}
+\eta_{ij}e^{f_i+f_j},\ \text{if}\ 1+\alpha_ie^{2f_i}<0,\ 1+\alpha_je^{2f_j}<0,
\end{gather*}
with
\begin{equation*}
\coth d_{ij}
=\frac{\alpha_ie^{2f_i}}{\sinh l_{ij}}
\sqrt{\frac{1+\alpha_je^{2f_j}}{1+\alpha_ie^{2f_i}}}
+\frac{\eta_{ij}e^{f_i+f_j}}{\sinh l_{ij}}.
\end{equation*}
This gives the discrete conformal structures (\ref{Eq: DCS4}) and (\ref{Eq: new 2}).
\end{description}
It is easy to check that the $\alpha_i$, $\eta_{ij}$ derived in any two adjacent right-angled hyperbolic  hexagons are equal.
\end{description}
Now we prove the claims (a), (b) and (c) in Theorem \ref{Thm: DCS}.
\begin{description}
\item[(1)]
By the arguments above, the discrete conformal structures (\ref{Eq: DCS3}), (\ref{Eq: new 1}) and (\ref{Eq: DCS1}) exist on $(\Sigma,\mathcal{T})$ only if $\frac{\sinh d_{ij}}{\sinh d_{ji}}>0$ for the edge $\{ij\}$.
Thus it is possible that there exist only the discrete conformal structures (\ref{Eq: DCS3}), (\ref{Eq: new 1}) or (\ref{Eq: DCS1}) on an ideally triangulated surface with boundary.
However, the discrete conformal structures (\ref{Eq: DCS3}), (\ref{Eq: new 1}) and (\ref{Eq: DCS1}) can not exist simultaneously in pairs on the same ideally triangulated surface with boundary.

Firstly, the different signs of $1+\alpha_re^{2f_r}$ ensure that (\ref{Eq: DCS3}) and (\ref{Eq: new 1}) can not exist simultaneously on the same ideally triangulated surface with boundary.
Otherwise, due to the connectedness of $(\Sigma,\mathcal{T})$,
there exists a right-angled hyperbolic hexagon $\{ijk\}\in F$ such that $l_{ij}$ is given by (\ref{Eq: DCS3}) and $l_{jk}$ is given by (\ref{Eq: new 1}).
Then $1+\alpha_je^{2f_j}>0$ and $1+\alpha_je^{2f_j}<0$.
This is a contradiction.

Secondly,
suppose that there exists a right-angled hyperbolic hexagon $\{ijk\}$ with $l_{ij}$ given by (\ref{Eq: DCS1}).
Without loss of generality, we can assume $\frac{\partial H}{\partial f_i}\equiv2,\
\frac{\partial H}{\partial f_j}\equiv-2$, otherwise change the orientation of the edge $\{ij\}$.
By the arguments below (\ref{Eq: F6}),
we see that $l_{jk},l_{ik}$ are also given by (\ref{Eq: DCS1}).
Similar arguments also apply to the discrete conformal structures (\ref{Eq: DCS3}) and (\ref{Eq: new 1}).
Hence, (\ref{Eq: DCS3}) and (\ref{Eq: DCS1}), or (\ref{Eq: new 1}) and (\ref{Eq: DCS1}) can not exist simultaneously on the same ideally triangulated surface with boundary.
This gives the proof of the claim (a).

\item[(2)]
By the arguments above again, the discrete conformal structures (\ref{Eq: DCS4}), (\ref{Eq: new 2}) and (\ref{Eq: DCS2}) exist in $(\Sigma,\mathcal{T})$ only if $\frac{\sinh d_{ij}}{\sinh d_{ji}}<0$  for the edge $\{ij\}$.
This implies that the discrete conformal structures (\ref{Eq: DCS4}), (\ref{Eq: new 2}) and (\ref{Eq: DCS2}) can not exist alone on an ideally triangulated surface with boundary.
Otherwise, we have $\frac{\sinh d_{ij}}{\sinh d_{ji}}<0$,\ $\frac{\sinh d_{ki}}{\sinh d_{ik}}<0$ and $\frac{\sinh d_{jk}}{\sinh d_{kj}}<0$
for a right-angled hyperbolic hexagon $\{ijk\}$, which implies $\frac{\sinh d_{ij}\sinh d_{jk}\sinh d_{ki}}{\sinh d_{ji}\sinh d_{ik}\sinh d_{kj}}<0$.
This contradicts the condition (\ref{Eq: compatible condition}).
Moreover, the discrete conformal structures (\ref{Eq: DCS4}), (\ref{Eq: new 2}) and (\ref{Eq: DCS2}) can not exist simultaneously in pairs on the same ideally triangulated surface with boundary.
The proof is the same as that for the second part of claim (a).
This gives the proof of the claim (b).

\item[(3)]
By the condition (\ref{Eq: compatible condition}), if one of $\frac{\sinh d_{ij}}{\sinh d_{ji}}$,\ $\frac{\sinh d_{ki}}{\sinh d_{ik}}$,  $\frac{\sinh d_{jk}}{\sinh d_{kj}}$ is negative,
then there is another unique one of $\frac{\sinh d_{ij}}{\sinh d_{ji}}$,\ $\frac{\sinh d_{ki}}{\sinh d_{ik}}$, $\frac{\sinh d_{jk}}{\sinh d_{kj}}$ is also negative.
For example,  $\frac{\sinh d_{ij}}{\sinh d_{ji}}<0$, $\frac{\sinh d_{ki}}{\sinh d_{ik}}<0$ and $\frac{\sinh d_{jk}}{\sinh d_{kj}}>0$.
This implies that two of the edge lengths $l_{ij}, l_{jk}, l_{ki}$ of a right-angled hyperbolic hexagon $\{ijk\}$ are induced by (\ref{Eq: DCS4}), (\ref{Eq: new 2}) or (\ref{Eq: DCS2}),
and the last one is induced by (\ref{Eq: DCS3}), (\ref{Eq: new 1}) or (\ref{Eq: DCS1}).
By the claim (b), the two of $l_{ij}, l_{jk}, l_{ki}$ with negative $\frac{\sinh d_{st}}{\sinh d_{ts}}$ are induced by the same one of (\ref{Eq: DCS4}), (\ref{Eq: new 2}) and (\ref{Eq: DCS2}).
By permutations and combinations, there are nine mixed types of discrete conformal structures, i.e.,
the discrete conformal structures (\ref{Eq: DCS3}) (\ref{Eq: new 2}), (\ref{Eq: DCS3}) (\ref{Eq: DCS2}), (\ref{Eq: new 1}) (\ref{Eq: DCS4}), (\ref{Eq: new 1}) (\ref{Eq: DCS2}), (\ref{Eq: DCS1}) (\ref{Eq: DCS4}), (\ref{Eq: DCS1}) (\ref{Eq: new 2}), (\ref{Eq: DCS3}) (\ref{Eq: DCS4}), (\ref{Eq: new 1}) (\ref{Eq: new 2}) and (\ref{Eq: DCS1}) (\ref{Eq: DCS2}).
The first six mixed types of discrete conformal structures can not exist on the ideally triangulated surface with boundary $(\Sigma,\mathcal{T})$,
while the last three can.
The proof is also the same as that for the second part of claim (a).
This proves the claim (c).
\end{description}
\qed

\begin{remark}
The trick of the proof of Theorem \ref{Thm: DCS} is taken from Glickenstein-Thomas \cite{GT}.
By similar calculations, we give a full classification of the discrete conformal structures on closed surfaces in \cite{X-Z}, which completes Glickenstein-Thomas's classification of  the discrete conformal structures on closed surfaces in \cite{GT}.
\end{remark}

\section{Relationships with 3-dimensional hyperbolic geometry}\label{section 6}

The relationships between discrete conformal structures on closed surfaces and 3-dimensional hyperbolic geometry were first observed by Bobenko-Pinkall-Springborn \cite{BPS} in the case of Luo’s vertex scaling of piecewise linear metrics.
Motivated by Bobenko–Pinkall–Springborn’s observations \cite{BPS}, Zhang-Guo-Zeng-Luo-Yau-Gu \cite{Zhang-Guo-Zeng-Luo-Yau-Gu} further constructed 18 types of discrete conformal structures on closed surfaces in different background geometries by perturbing the ideal vertices of an ideal hyperbolic tetrahedron to be hyperbolic, ideal or hyper-ideal.
Recently, Xu \cite{Xu22} constructed a new class of discrete conformal structures on surfaces with boundary by perturbing the ideal vertices of the ideal hyperbolic tetrahedron to be hyper-ideal such that some edges of the resulting hyperbolic tetrahedron do not intersect with $\mathbb{H}^3\cup \partial \mathbb{H}^3$.
Motivated by the above works, we construct ten types of generalized hyperbolic tetrahedra and
give the relationships between the discrete conformal structures on surfaces with boundary in Theorem \ref{Thm: DCS} and the 3-dimensional hyperbolic geometry.

\begin{figure}[!ht]
\centering
\includegraphics[scale=0.8]{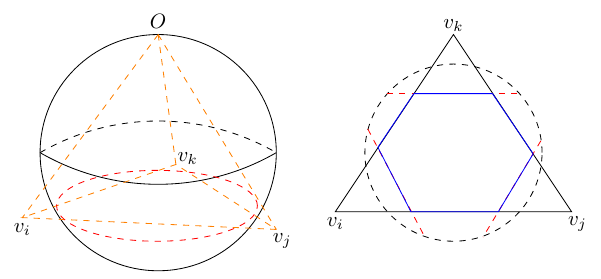}
\caption{A generalized hyperbolic tetrahedron $Ov_iv_jv_k$ with an ideal vertex $O$ and three hyper-ideal vertices $v_i,v_j,v_k$, and its bottom hyper-ideal hyperbolic triangle $v_iv_jv_k$ in the Kelin model.}
\label{Figure_7}
\end{figure}

Suppose $Ov_iv_jv_k$ is an ideal hyperbolic tetrahedron in $\mathbb{H}^3$ with four ideal vertices $O,v_i,v_j,v_k$.
Each ideal vertex is attached with a horosphere.
We perturb the ideal vertices $v_i,v_j,v_k$ to be hyper-ideal such that the line segments $v_iv_j, v_jv_k, v_iv_k$ intersect with the 3-dimensional hyperbolic space $\mathbb{H}^3$ in the Klein model.
Then the intersection of the hyper-ideal hyperbolic triangle $v_iv_jv_k$ with $\mathbb{H}^3$ corresponds to a right-angled hyperbolic hexagon.
Please refer to Figure \ref{Figure_7}.

\begin{figure}[!ht]
\centering
\includegraphics[scale=1]{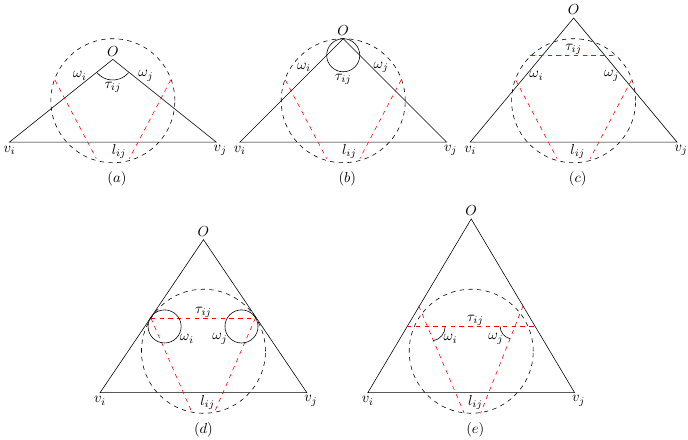}
\caption{Lateral generalized triangles in the Kelin model.}
\label{Figure_6}
\end{figure}

We first give the geometric explanations of the discrete conformal structures (\ref{Eq: DCS3}), (\ref{Eq: new 1}) and (\ref{Eq: DCS1}).
To make our statement more clear, we use a lateral generalized triangle $Ov_iv_j$ of the generalized hyperbolic tetrahedron $Ov_iv_jv_k$.
Please refer to Figure \ref{Figure_6}.
Denote the hyperbolic line dual to $O,v_i,v_j$ by $L_O,L_i,L_j$ respectively.
According to Theorem \ref{Thm: DCS} (a), the discrete conformal structures (\ref{Eq: DCS3}), (\ref{Eq: new 1}) and (\ref{Eq: DCS1}) can exist alone on an ideally triangulated surface with boundary.
Hence, the three lateral generalized triangles with the cosine law from (\ref{Eq: DCS3}), (\ref{Eq: new 1}) or (\ref{Eq: DCS1}) of a hyperbolic tetrahedron are the same.

Keeping the hyper-ideal vertices $v_i,v_j,v_k$ fixed, we start to perturb the ideal vertex $O$.

\textbf{(I):}
We first perturb the ideal vertex $O$ to be hyperbolic, then the triangle $Ov_iv_j$ forms a generalized hyperbolic triangle with two hyper-ideal vertices and one hyperbolic vertex.
Please refer to Figure \ref{Figure_6} ($a$).
Denote the angle at $O$ in the triangle $Ov_iv_j$ by $\tau_{ij}$.
There exists a unique hyperbolic segment $L_{ij}$ perpendicular to $L_i$ and $L_j$, the length of which is denoted by $l_{ij}$.
The distance from $O$ to $L_i$ and $L_j$ is denoted by $\omega_i$ and $\omega_j$ respectively.
By $(-1,-1,1)$ type generalized hyperbolic cosine law in Guo-Luo \cite{GL2}, we have
\begin{equation*}
\cosh l_{ij}=\sinh \omega_i\sinh \omega_j-\cos\tau_{ij}\cosh \omega_i\cosh \omega_j.
\end{equation*}
This corresponds to the discrete conformal structure (\ref{Eq: new 1}) by letting $e^f=\cosh\omega$ and $\eta=-\cos\tau$.

\textbf{(II):}
We perturb the ideal vertex $O$ to be ideal such that the line segments $Ov_i, Ov_j, Ov_k$ intersect with the 3-dimensional hyperbolic space $\mathbb{H}^3$.
Then the generalized triangle $Ov_iv_j$ is a generalized hyperbolic triangle with two hyper-ideal vertices and one ideal vertex.
Please refer to Figure \ref{Figure_6} ($b$).
The length of the intersection of the horocycle at $O$ bounded by $Ov_i$ and $Ov_j$ is denoted by $\tau_{ij}$.
The distance from the horocycle at $O$ to $L_i$ and $L_j$ is denoted by $\omega_i$ and $\omega_j$ respectively.
By the $(-1,-1,0)$ type generalized hyperbolic cosine law in Guo-Luo \cite{GL2}, we have
\begin{equation*}
\cosh l_{ij}
=-\cosh(\omega_i-\omega_j)
+\frac{1}{2}\tau^2_{ij}e^{\omega_i+\omega_j}.
\end{equation*}
This corresponds to the discrete conformal structure (\ref{Eq: DCS1}) with $C_{ij}=0$ by letting $f=\omega$ and $\eta=\frac{1}{2}\tau^2$.

\textbf{(III):} We continue to perturb the ideal vertex $O$ to be hyper-ideal such that the line segments $Ov_i, Ov_j, Ov_k$ intersect with the 3-dimensional hyperbolic space $\mathbb{H}^3$.
The intersection of the hyperbolic plane $P_O$ dual to $O$ with $Ov_iv_jv_k$ is a hyperbolic triangle in $\mathbb{H}^3$.
The lines $L_O,L_i,L_j$ form a right-angled hyperbolic  hexagon. Please refer to Figure \ref{Figure_6} ($c$).
Furthermore, there exists a unique hyperbolic segment $L_{Oi}$ perpendicular to $L_O$ and $L_i$, the length of which is denoted by $\omega_{i}$.
And there exists a unique hyperbolic segment $L_{Oj}$ perpendicular to $L_O$ and $L_j$, the length of which is denoted by $\omega_{j}$.
And the length of $L_O$ between  $L_{Oi}$ and $L_{Oj}$ is denoted by $\tau_{ij}$.
By the $(-1,-1,-1)$ type generalized hyperbolic cosine law in Guo-Luo \cite{GL2}, we have
\begin{equation*}
\cosh l_{ij}=-\cosh \omega_i\cosh \omega_j+\cosh\tau_{ij}\sinh \omega_i\sinh \omega_j.
\end{equation*}
This corresponds to the discrete conformal structure (\ref{Eq: DCS3}) by letting $\alpha\equiv1$, $e^{f}=\sinh \omega$ and $\eta=\cosh\tau$.

\textbf{(IV):}  We further perturb the hyperbolic vertex $O$ such that the line segments $Ov_i, Ov_j, Ov_k$ are tangential to $\partial\mathbb{H}^3$,
then the intersection of the hyperbolic plane $P_O$ dual to $O$ with $Ov_iv_jv_k$ is an ideal hyperbolic triangle with three ideal vertices.
Since the line $L_O$ intersects the line $L_i$ and the line $L_j$ on $\partial\mathbb{H}^2$ respectively,
then the lines $L_O,L_i,L_j$ form a generalized hyperbolic triangle with one hyper-ideal vertex and two ideal vertices. Please refer to Figure \ref{Figure_6} ($d$).
Furthermore, the length of the intersection of the associated horocycle bounded by $L_i$ and $L_O$ is denoted by $\omega_i$,
and the length of the intersection of the associated horocycle bounded by $L_j$ and $L_O$ is denoted by $\omega_j$.
The distance of the two horocycles is denoted by $\tau_{ij}$.
By the $(0,0,-1)$ type generalized hyperbolic cosine law in Guo-Luo \cite{GL2}, we have
\begin{equation*}
\cosh l_{ij}=-1+2e^{\tau_{ij}}\omega_i\omega_j,
\end{equation*}
which corresponds to the discrete conformal structure (\ref{Eq: DCS3}) by letting $\alpha\equiv0$, $e^{f}=\omega$ and $\eta=2e^{\tau}$.
This is exactly Guo's vertex scalings on surfaces with boundary in \cite{Guo}.

\textbf{(V):} We perturb the ideal vertex $O$ such that the line segments $Ov_i, Ov_j, Ov_k$ do not intersect with $\mathbb{H}^3\cup\partial\mathbb{H}^3$,
then the intersection of the hyperbolic plane $P_O$ dual to $O$ with $Ov_iv_jv_k$ corresponds to a right-angled hyperbolic hexagon.
The lines $L_O,L_i,L_j$ can form a generalized hyperbolic triangle with one hyper-ideal vertex and two hyperbolic vertices. Please refer to Figure \ref{Figure_6} ($e$).
Furthermore, the intersection angles of $L_O$ and $L_i$, $L_j$ are denoted by $\omega_i$ and $\omega_j$ respectively.
The length of $L_O$ between $L_i$ and $L_j$ is denoted by $\tau_{ij}$.
By the $(1,1,-1)$ type generalized hyperbolic cosine law in Guo-Luo \cite{GL2}, we have
\begin{equation*}
\cosh l_{ij}=-\cos\omega_i\cos\omega_j+\cosh \tau_{ij}\sin\omega_i\sin\omega_j.
\end{equation*}
This corresponds to the discrete conformal structure (\ref{Eq: DCS3}) by letting $\alpha\equiv-1$, $e^{f}=\sin \omega$ and $\eta=\cosh \tau$.
This is one of Xu's discrete conformal structures on surfaces with boundary in \cite{Xu22}.

\begin{figure}[!ht]
\centering
\includegraphics[scale=1]{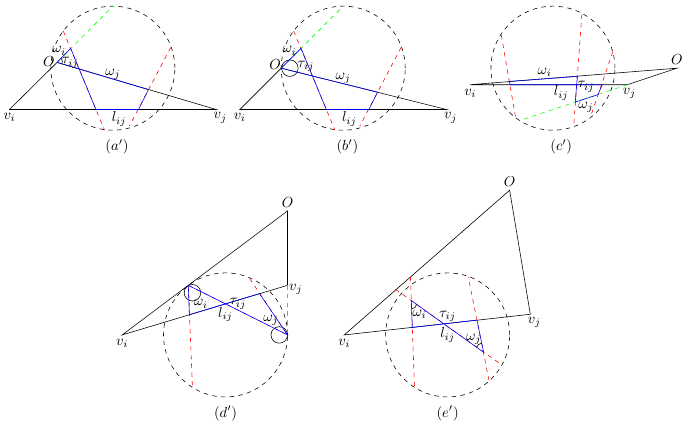}
\caption{Lateral generalized triangles in the Kelin model and the twisted triangles are bounded by blue lines.}
\label{Figure_9}
\end{figure}

Now we give the geometric explanations of the discrete conformal structures (\ref{Eq: DCS4}), (\ref{Eq: new 2}) and (\ref{Eq: DCS2}).
According to Theorem \ref{Thm: DCS} (b),
the discrete conformal structures (\ref{Eq: DCS4}), (\ref{Eq: new 2}) and (\ref{Eq: DCS2}) can not exist alone on an ideally triangulated surface with boundary.
And we consider the mixed type of discrete conformal structures, i.e., the discrete hyperbolic metrics $l_{ij}, l_{jk}, l_{ki}$ of a right-angled hyperbolic hexagon $\{ijk\}\in F$ are induced by (\ref{Eq: DCS3}) and (\ref{Eq: DCS4}), (\ref{Eq: new 1}) and (\ref{Eq: new 2}) or (\ref{Eq: DCS1}) and (\ref{Eq: DCS2}).
To ensure the condition (\ref{Eq: compatible condition}) holds,
there must be one of $l_{ij}, l_{jk}, l_{ki}$ from (\ref{Eq: DCS3}) and the other two from (\ref{Eq: DCS4}),
one from (\ref{Eq: new 1}) and the other two from (\ref{Eq: new 2}),
or one from (\ref{Eq: DCS1}) and the other two from (\ref{Eq: DCS2}).
In other words, in a hyperbolic tetrahedron,
there exist one lateral generalized triangle with the cosine law from (\ref{Eq: DCS3}) and the other two from (\ref{Eq: DCS4}),
one lateral generalized triangle with the cosine law from (\ref{Eq: new 1}) and the other two from (\ref{Eq: new 2}),
or one lateral generalized triangle with the cosine law from (\ref{Eq: DCS1}) and the other two from (\ref{Eq: DCS2}).
For the discrete conformal structures (\ref{Eq: DCS3}), (\ref{Eq: new 1}) and (\ref{Eq: DCS1}),
the lateral generalized triangles are still shown in Figure \ref{Figure_6}.
However, for the discrete conformal structures (\ref{Eq: DCS4}), (\ref{Eq: new 2}) and (\ref{Eq: DCS2}),
the lateral generalized triangles are some twisted hyperbolic triangles introduced by Roger-Yang \cite{Roger-Yang}.
Please refer to Figure \ref{Figure_9}.
One can also refer to \cite{B-Y,Roger-Yang} for more information on twisted hyperbolic triangles.

\begin{figure}[!ht]
\centering
\includegraphics[scale=0.9]{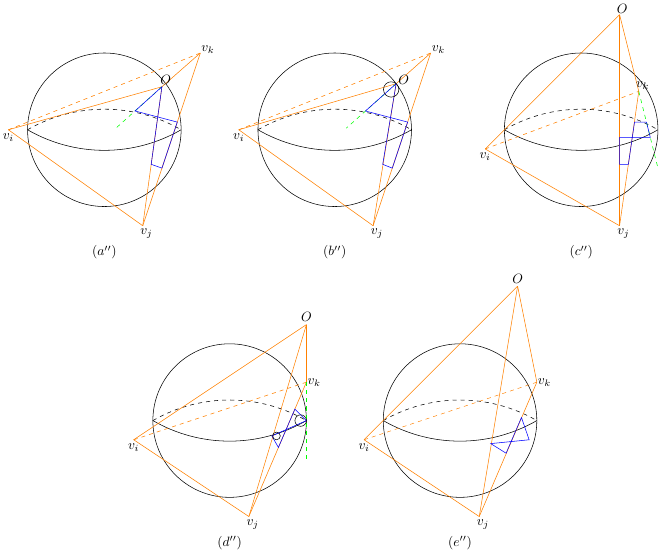}
\caption{Schematic diagram of four generalized hyperbolic tetrahedra and the corresponding lateral generalized triangles in the Kelin model.}
\label{Figure_10}
\end{figure}

\textbf{(VI):}
For the discrete conformal structures (\ref{Eq: new 2}),
the lateral generalized triangle is shown in Figure \ref{Figure_9} ($a^\prime$).
Unlike the case where $L_i$ intersect the line segment $Ov_i$ as shown in Figure \ref{Figure_6} (a),
in this case, $L_i$ intersects the reverse extension (denoted by $v_iO$) of the line segment $Ov_i$.
There exists a unique hyperbolic segment $L_{ij}$ perpendicular to $L_i$ and $L_j$, the length of which is denoted by $l_{ij}$.
The angle at $O$ bounded by $v_iO$ and $Ov_j$ is denoted by $\tau_{ij}$.
The distance from $O$ to $L_i$ and $L_j$ is denoted by $\omega_i$ and $\omega_j$ respectively.
By the cosine laws of $(1,-1,-1)$ type twisted generalized triangle in Roger-Yang \cite{Roger-Yang},
we have
\begin{equation*}
\cosh l_{ij}
=-\sinh\omega_i\sinh\omega_j
+\cos\tau_{ij}\cosh\omega_i\cosh\omega_j.
\end{equation*}
This corresponds to the discrete conformal structure (\ref{Eq: new 2}) by letting $e^f=\cosh\omega$ and $\eta=\cos\tau$.
The hyperbolic tetrahedron $Ov_iv_jv_k$ with one lateral generalized triangles shown in Figure \ref{Figure_6} ($a$) and two shown in Figure \ref{Figure_9} ($a^\prime$) is depicted in Figure \ref{Figure_10} ($a^{\prime\prime}$).
Specially, the generalized triangle $Ov_iv_j$ is shown in Figure \ref{Figure_6} ($a$), and $Ov_jv_k$ and $Ov_iv_k$ are shown in Figure \ref{Figure_9} ($a^\prime$).

\textbf{(VII):} For the discrete conformal structures (\ref{Eq: DCS2}),
the lateral generalized triangle is shown in Figure \ref{Figure_9} ($b^\prime$).
Unlike the case where the line segment $Ov_i$ intersects $\mathbb{H}^3$ as shown in Figure \ref{Figure_6} (b),
in this case, the reverse extension (denoted by $v_iO$) of the line segment $Ov_i$ intersects $\mathbb{H}^3$.
The length of the intersection of the horocycle at $O$ bounded by $v_iO$ and $Ov_j$ is denoted by $\tau_{ij}$.
The distance from the horocycle at $O$ to $L_i$ and $L_j$ is denoted by $\omega_i$ and $\omega_j$ respectively.
By the cosine laws of $(0,-1,-1)$ type twisted generalized triangle in Roger-Yang \cite{Roger-Yang},
we have
\begin{equation*}
\cosh l_{ij}
=\cosh(\omega_i-\omega_j)
-\frac{1}{2}\tau^2_{ij}e^{\omega_i+\omega_j}.
\end{equation*}
This corresponds to the discrete conformal structure (\ref{Eq: DCS2}) by letting $f=\omega$ and $\eta=-\frac{1}{2}\tau^2$.
The hyperbolic tetrahedron $Ov_iv_jv_k$ with one lateral generalized triangles shown in Figure \ref{Figure_6} ($b$) and two shown in Figure \ref{Figure_9} ($b^\prime$) is depicted in Figure \ref{Figure_10} ($b^{\prime\prime}$).
Specially, the generalized triangle $Ov_iv_j$ is shown in Figure  \ref{Figure_6} ($b$), and $Ov_jv_k$ and $Ov_iv_k$ are shown in Figure \ref{Figure_9} ($b^\prime$).

\textbf{(VIII):}
Figure \ref{Figure_9} ($c^\prime$) corresponds to the discrete conformal structure (\ref{Eq: DCS4}) for $\alpha\equiv1$, $e^{f}=\sinh \omega$ and $\eta=\cosh\tau$.
Unlike the case where the line segment $Ov_j$ intersects $\mathbb{H}^3$ as shown in Figure \ref{Figure_6} (c),
in this case, the extension (still denoted by $Ov_j$) of the line segment $Ov_j$ intersects $\mathbb{H}^3$.
There exists a unique hyperbolic segment $L_{Oi}$ perpendicular to $L_O$ and $L_i$, the length of which is denoted by $\omega_{i}$.
And there exists a unique hyperbolic segment $L_{Oj}$ perpendicular to $L_O$ and $L_j$, the length of which is denoted by $\omega_{j}$.
The length of $L_O$ between $L_{Oi}$ and $L_{Oj}$ is denoted by $\tau_{ij}$.
By the cosine laws of $(-1,-1,-1)$ type twisted generalized triangle in Roger-Yang \cite{Roger-Yang},
we have
\begin{equation*}
\cosh l_{ij}=\cosh \omega_i\cosh \omega_j+\cosh \tau_{ij}\sinh \omega_i\sinh \omega_j.
\end{equation*}
This corresponds to the discrete conformal structure (\ref{Eq: DCS4}) by letting $\alpha\equiv1$, $e^{f}=\sinh \omega$ and $\eta=\cosh\tau$.
The hyperbolic tetrahedron $Ov_iv_jv_k$ with one lateral generalized triangles shown in Figure  \ref{Figure_6} ($c$) and two shown in Figure \ref{Figure_9} ($c^\prime$) is depicted in Figure \ref{Figure_10} ($c^{\prime\prime}$).
Specially, the generalized triangle $Ov_iv_j$ is shown in Figure  \ref{Figure_6} ($c$), and $Ov_jv_k$ and $Ov_iv_k$ are shown in Figure \ref{Figure_9} ($c^\prime$).
One can also refer to Figure 10 ($b$) in Belletti-Yang \cite{B-Y} for generalized hyperbolic tetrahedron $Ov_iv_jv_k$,
where $v_l=v_i,v_i=v_j,v_j=O$.

\textbf{(IX):}
Figure \ref{Figure_9} ($d^\prime$) corresponds to the discrete conformal structure (\ref{Eq: DCS4}) for $\alpha\equiv0$, $e^{f}=\omega$ and $\eta=2e^{\tau}$.
Unlike the case where the line segment $Ov_j$ is tangential to $\partial\mathbb{H}^3$ as shown in Figure \ref{Figure_6} (d),
in this case, the extension (denoted by $Ov_j$) of the line segment $Ov_j$ is tangential to $\partial\mathbb{H}^3$.
The length of the intersection of the associated horocycle bounded by $L_i$ and $L_O$ is denoted by $\omega_i$,
and the length of the intersection of the associated horocycle bounded by $L_j$ and $L_O$ is denoted by $\omega_j$.
The distance of the two horocycles is denoted by $\tau_{ij}$.
By the cosine laws of $(0,0,-1)$ type twisted generalized triangle in Roger-Yang \cite{Roger-Yang},
we have
\begin{equation*}
\cosh l_{ij}=1+2e^{\tau_{ij}}\omega_i\omega_j.
\end{equation*}
This corresponds to the discrete conformal structure (\ref{Eq: DCS4}) by letting $\alpha\equiv0$, $e^{f}=\omega$ and $\eta_{ij}=2e^{\tau_{ij}}$.
The hyperbolic tetrahedron $Ov_iv_jv_k$ with one lateral generalized triangles shown in Figure  \ref{Figure_6} ($d$) and two shown in Figure \ref{Figure_9} ($d^\prime$) is depicted in Figure \ref{Figure_10} ($d^{\prime\prime}$).
Specially, the generalized triangle $Ov_iv_j$ is shown in Figure  \ref{Figure_6} ($d$), and $Ov_jv_k$ and $Ov_iv_k$ are shown in Figure \ref{Figure_9} ($d^\prime$).

\textbf{(X):}
Figure \ref{Figure_9} ($e^\prime$) corresponds to the discrete conformal structure (\ref{Eq: DCS4}) for $\alpha\equiv-1$, $e^{f}=\sin \omega$ and $\eta=\cosh \tau$.
Unlike the case where $L_O$ does not intersect $L_{ij}$ as shown in Figure \ref{Figure_6} (e),
in this case, $L_O$ intersects $L_{ij}$.
The intersection angles of $L_O$ and $L_i$, $L_j$ are denoted by $\omega_i$ and $\omega_j$ respectively.
The length of $L_O$ between $L_i$ and $L_j$ is denoted by $\tau_{ij}$.
By the cosine laws of $(1,1,-1)$ type twisted generalized triangle in Roger-Yang \cite{Roger-Yang},
we have
\begin{equation*}
\cosh l_{ij}=\cos \omega_i\cos \omega_j+\cosh \tau_{ij}\sin \omega_i\sin \omega_j.
\end{equation*}
This corresponds to the discrete conformal structure (\ref{Eq: DCS4}) by letting $\alpha\equiv -1$, $e^{f}=\sin\omega$ and $\eta=\cosh \tau$. The hyperbolic tetrahedron $Ov_iv_jv_k$ with one lateral generalized triangles shown in Figure  \ref{Figure_6} ($e$) and two shown in Figure \ref{Figure_9} ($e^\prime$) is depicted in Figure \ref{Figure_10} ($e^{\prime\prime}$).
Specially, the generalized triangle $Ov_iv_j$ is shown in Figure  \ref{Figure_6} ($e$), and $Ov_jv_k$ and $Ov_iv_k$ are shown in Figure \ref{Figure_9} ($e^\prime$).

\end{document}